\documentclass[12pt,a4paper,twoside,final,notitlepage, leqno]{article}
\usepackage[english]{babel}
\usepackage[T1]{fontenc}  
\usepackage{graphicx}
\usepackage{amsmath,amsthm,epsfig,amsfonts,bbm}  

\def \smb {{\scriptstyle \bullet }}
\newcommand{\monitem}{ \smallskip \noindent $\bullet$ \quad  } 
\newcommand{\moneq}{\vspace*{-6pt} \begin{equation} \displaystyle } 
\newcommand{\moneqstar}{\vspace*{-6pt} \begin{equation*} \displaystyle } 
\newcommand{\monendstar}{\vspace*{-6pt} \end{equation*}   }
\newcommand{\monend}{\vspace*{-6pt} \end{equation}   }

\def\R{{\rm I}\! {\rm R}}

\def\section*#1{}

\def\resume{\if@twocolumn
\section*{R\'esum\'e}
\else \small
\quotation{\bf \it R\'esum\'e \rule[1mm]{1.5mm}{0.2mm}\vspace{0pt}}
\fi}
\def\endresume{\if@twocolumn\else\endquotation\fi}
\def\abstract{\if@twocolumn
\noindent\section*{{\bf Abstract}}
\else \small
\quotation{\noindent \bf {Abstract.} \rule[1mm]{1.5mm}{0.2mm}\vspace{0pt}}
\fi}
\def\endabstract{\if@twocolumn\else\endquotation\fi}

\usepackage{fancyhdr}
\renewcommand{\headrulewidth}{0pt}
\begin{document}

\fancypagestyle{plain}{ \fancyfoot{} \renewcommand{\footrulewidth}{0pt}}
\fancypagestyle{plain}{ \fancyhead{} \renewcommand{\headrulewidth}{0pt}}

\bibliographystyle{alpha}


\title{\bf \LARGE  On rotational invariance \\ 
of lattice Boltzmann schemes
}

\author { {\large   Adeline Augier$^1$, Fran\c{c}ois Dubois$^{1,}$$^2$$^*$,     
Benjamin Graille$^1$,  Pierre Lallemand$^{3}$ } \\ ~ \\   
{\it \small $^1$    Department of Mathematics, University  Paris Sud, Orsay, France.}\\ 
{\it \small $^2$   Conservatoire National des Arts et M\'etiers, Paris, France,  } \\ 
{\it \small     Structural Mechanics and Coupled Systems Laboratory.  } \\   
{\it \small $^3$   Beiging Science Computing Research Center, China.}\\ 
{ \rm  \small  $^*$ corresponding author: francois.dubois@math.u-psud.fr . } ~ \\  }

\date  {~ \\ ~ \\  {  \rm  30  June 2013} 
 \footnote {\rm  \small $\,\,$ Contribution available online 9 July 2013  in 
 {\it Computers and Mathematics with Applications}, 
doi: http://dx.doi.org/10.1016/j.camwa.2013.06009.   
Invited communication presented at the  
 Ninth International Conference for Mesoscopic   Methods in Engineering 
and Science,  Taipei, Taiwan,  26 July 2012. Edition 02 January 2014.  }}

\maketitle
\renewcommand{\baselinestretch}{1.}

  \vspace{-.5cm}  
  
 \bigskip \bigskip  
\noindent {\bf  \large   Abstract  } 
 
\noindent  
We propose the derivation of acoustic-type isotropic partial 
differential equations that are equivalent
to linear lattice Boltzmann schemes with a density scalar   field and a momentum 
vector field as conserved moments. 
The corresponding linear equivalent  partial differential equations
are generated with a new ``Berliner version'' of the Taylor
expansion method. The  details  of the 
implementation are  presented. 
These ideas are applied for the D2Q9, D2Q13, 
D3Q19 and D3Q27    lattice Boltzmann schemes. 
Some limitations associated with necessary stability conditions 
are also presented.

 $ $ \\[.4mm]

 $ $ \\[.4mm]

\newpage 
\fancyfoot[C]{\oldstylenums{\thepage}}

~

 \noindent {\bf \large 1) \quad    Introduction  } 

 \monitem     
  Partial differential equations  like Navier-Stokes equations  
are    invariant by rotation  
 and all  space directions are equivalent.  
Due to the use of a given mesh, 
lattice Boltzmann schemes
cannot be  completely invariant by rotation.             
This difficulty was present in the early ages of lattice 
 gas automata. The initial model of Hardy, de Pazzis and Pomeau
\cite{HPP73} 
proposed  very impressive qualitative results but the associated   
fluid tensor was not  invariant by rotation. 
With a triangular mesh, the  second model of 
Frisch, Hasslacher and  Pomeau  \cite{FHP86} gives the correct physics.   
The lattice Boltzmann scheme with multiple   relaxation times   
is the fruit of the work of 
    Higuera and  Jim\'enez \cite{HJ89}, 
Higuera, Succi and Benzi \cite{HSB89},  
Qian, d'Humi\`eres and Lallemand \cite{QHL92}  
and  d'Humi\`eres \cite {DdH92}. 
It  uses in general square 
meshes and leads to isotropic physics for second order equivalent model   
as  analyzed in \cite{LL00}.   
The question of rotational invariance is still present in the lattice
Boltzmann community and  
a detailed analysis of moment isotropy has been proposed by 
Chen and  Orszag \cite {CH11}.  

 \monitem     
The invariance by rotation has to be kept as  much   
as possible in order to respect the     
correct propagation of waves. 
In \cite {DL11}, using the Taylor expansion method proposed     
in \cite {DL09} for general applications, 
we have developed a methodology to enforce a lattice Boltzmann scheme to
simulate correctly the physical acoustic waves up to fourth order of accuracy.
But unfortunately, stability is in general  guaranteed only 
if the viscosities of the waves are much  higher than authorized        
by common physics. In this contribution, we relax this  constraint and suppose that 
the equivalent partial differential equation of the scheme are invariant
under rotations. 
The objective of this contribution is to propose a methodology     
to fix the  parameters of lattice Boltzmann schemes                   
in order to ensure the invariance by rotation at a given order.

 \monitem     
  The outline is the following.
In Section~2, we consider the question of the 
algebraic form of linear high order  ``acoustic-type''
partial differential equations that are invariant by two dimensional and      
three-dimensional rotations. 
In Section~3, we recall the essential  properties  concerning multiple         
  relaxation times lattice Boltzmann schemes. 
The equivalent equation of a lattice Boltzmann scheme introduces 
naturally the notion of rotational invariance
at a given order. 
In the following sections, we develop a methodology to force acoustic-type 
lattice Boltzmann models to be invariant under rotations. We consider the
 four lattice Boltzmann schemes D2Q9, D2Q13, D3Q19 and D3Q27 in the sections~4
to~7. A conclusion ends our contribution. 
The Appendix~1  presents  with details the 
implementation of the ``Berliner version'' of our algorithm      
to derive explicitely the equivalent partial differential equations. Some 
long formulas associated of specific results 
for D2Q13  and D3Q27 lattice Boltzmann schemes
are presented in Appendix~2.

  \bigskip \bigskip \newpage   \noindent {\bf \large 2) \quad  Invariance by rotation of
  acoustic-type equations   } 
 
 \monitem      
With the help of group theory, and we refer the reader  {\it e.g.} to 
Hermann Weyl \cite{We36} or  Goodman and  Wallach \cite{GW98},   
it is possible to write {\it a priori} the general form
of  systems of linear partial differential equations invariant by rotation. 
More precisely, if the unknown is composed by one scalar field $ \, \rho \,$ 
(invariant function under a rotation of the space) 
and one vector field $ \, J \,$ (a vector valued function that is transformed 
in a similar way to how   
  than cartesian coordinates are  
transformed when a rotation is applied),   
a linear  partial differential equation 
invariant by rotation  
is constrained in a strong manner. 
Using some fundamental aspects of group theory and in particular the Schur lemma
 (see  {\it e.g.}   Goodman and  Wallach \cite{GW98}), 
it is possible to prove that  general linear partial differential equations 
of acoustic type  that are invariant by rotation  admit  the   form described below. 

 \monitem      
In the bidimentional  case, we introduce the notation   
\moneq   \label{rot-perp-2d}   
\nabla^\perp \, \equiv \, \Big( {{\partial}\over{\partial x}} \,,\, 
{{\partial}\over{\partial y}} \Big)^\perp 
 \, =  \,  \Big( {{\partial}\over{\partial y}} \,,\, 
 - {{\partial}\over{\partial x}} \Big)  \,,  \quad  
J^\perp   \, \equiv \,   ( J_x , \, J_y)^\perp \, =  \, ( J_y, \, -J_x) \,.  \quad  
\monend  
Then acoustic type partial differential equations are of the form
%
%
\moneq   \label{rot-eqs-2d}   
\left\{\begin{array}{l}   \displaystyle 
\partial_t \rho + \sum_{k\geq0} \Big( \alpha_k \, \triangle^k\rho
+ \beta_k \, \triangle^k  {\rm div}  J + \gamma_k \, \triangle^k   {\rm div}  (J^\perp)
\Big) = 0 \,, 
 \\ \displaystyle   \vspace{-.5cm}  ~  \\ \displaystyle  
\displaystyle \partial_t J  +\sum_{k\geq0} \Big(  \delta_k \, \nabla\triangle^k\rho 
+ \mu_k  \triangle^k J + 
\zeta_k \,  \nabla {\rm div} \triangle^k J   
 \\ \displaystyle   \vspace{-.9cm}  ~  \\ \displaystyle
\qquad \qquad \qquad  \qquad  \qquad \qquad  
+ \, \varepsilon_k\, \nabla^\perp \triangle^k\rho + \, \nu_k \, \triangle^k J^\perp + 
\eta_k \, \nabla {\rm div} \triangle^kJ^\perp \Big) = 0 \, ,   
\end{array} \right. \monend  
where the real coefficients $ \, \alpha_k$,  $ \, \beta_k$,   $ \, \gamma_k $,  
$ \, \delta_k$,  $ \, \mu_k  $,   $ \,  \zeta_k $,  $ \,  \varepsilon_k $, 
$ \, \nu_k \,$ and  $ \, \eta_k \, $ are in finite number. 
The tridimensional case is essentially analogous. The ``acoustic type'' linear partial
differential equations invariant by rotation take necessarily the form 
%
%
\moneq   \label{rot-eqs-3d}      
\left\{\begin{array}{l}   \displaystyle  
\partial_t \rho +\sum_{k\geq0}  \Big( \alpha_k \, \triangle^k \rho 
+ \beta_k  \,  {\rm div} \triangle^k J  \Big) = 0 
 \\ \displaystyle   \vspace{-.5cm}  ~  \\ \displaystyle  
\displaystyle \partial_t J +\sum_{k\geq0} \Big( \delta_k \, \nabla \triangle^k \rho 
+ \mu_k  \,  \triangle^k J
+ \eta_k \,  \nabla {\rm div}  \triangle^k J   
+ \varphi_k \, \mbox{curl} \, \triangle^k J  \Big) \, = \, 0  \, ,   
\end{array} \right. \monend  
with an analogous convention  that the sums in the relations (\ref{rot-eqs-3d})  
contain only  a finite number of such terms. 

\bigskip \bigskip   \noindent {\bf \large 3) \quad  Lattice Boltzmann
  schemes  with multiple relaxation times }

 \monitem      
Each iteration of a lattice Boltzmann scheme is composed by two steps:  relaxation and
propagation. 
The relaxation is local in space: the particle distribution $ \, f(x) \in \R^q  \,$
for  $ \, x \,$ a node of the lattice $ \,  {\cal L} , \,$ is transformed 
into a ``relaxed''' distribution $ \, f^*(x)
\,$ that is  non linear in general.  
In this contribution, we restrict to linear functions 
$\, \R^q \ni f \longmapsto f^* \in \R^q  . \,$
As usual with the d'Humi\`eres scheme \cite{DdH92}, we introduce an invertible matrix $M$ with $q$
lines and $q$ columns. The  moments $m$ are obtained from the particle distribution thanks to the
associated transformation
\moneq   \label{moments}      
m_k \,=\, \sum_{j = 0}^{q-1} M_{k \, j} \, f_j \,, \qquad 0 \leq k \leq q-1 \, . 
\monend  
Then we consider the conserved moments $ \, W \in \R^N :$
\moneq   \label{conserves}
W_i \,=\, m_i \,, \qquad 0 \leq i \leq N-1 \, . 
\monend  
For the 
usual acoustic equations for $d$ space dimensions,  we have $ \, N = d+1 . \,$ 
The first moment is the density and the next ones are    
composed by the $d$ components of the physical momentum. 
Then we define a conserved value $\, m_k^{\rm eq} \,$ 
for the non-equilibrium moments $\, m_k \,$
for $ \, k \geq N .\,$   
With the help of ``Gaussian''  functions $ \, G_k(\smb), $     
we obtain:              
\moneq   \label{non-conserves}
m_k^{\rm eq}  \,=\, G_k (W) \,,  \qquad  N \leq k  \leq  q-1 \, . 
\monend  
In the present contribution, we suppose that this equilibrium value is a linear function
of the conserved variables. In other terms, the Gaussian functions are linear: 
\moneq   \label{gaussian}
G_{N+\ell} (W)   \,=\, \sum_{i=1}^{n-1} \, E_{\ell i} \, W_i \,, \qquad \ell \geq 0 \,  
\monend  
for some equilibrium coefficients $\, E_{\ell i} \,$ for $\, \ell \geq 0 \,$ and $ \, 0
\leq i \leq N-1 $.

 \monitem      
The relaxed moments $ \, m_k^* \,$ are  linear functions of $\, m_k \,$ and 
$ \, m_k^{\rm eq} :$ 
\moneq   \label{relaxation}
m_k^* \,= m_k + s_k \, \big( m_k^{\rm eq} - m_k \big) \,, \qquad k \geq N  \, .
\monend  
%
For a  stable  scheme, we   have      
\moneq   \label{stabi}
0 \, < \, s_k \, < \, 2 \, . 
\monend  
We remark that if $\, s_k=0 $, the corresponding moment is conserved. 
In some particular cases, the value $\, s_k =  2 \,$ can also be used
(see {\it e.g.} \cite{Du12, Du13}).        
The conserved moments are not affected by the relaxation: 
\moneqstar 
m_i^* \,=\, m_i \, = \, W_i \,, \qquad 0 \leq i \leq N-1 \, . 
\monendstar  
From the moments $ \, m_\ell^*  \,$ for $ \, 0 \leq \ell \leq q-1 \,$ we deduce the
particle distribution 
$ \, f_j^* \,$ by resolution of the linear system 
\moneqstar 
M \, f^* \,=\, m^* \, .  
\monendstar   

 \monitem   
The propagation step couples the node $\, x \in {\cal L} \,$ with his neighbours 
$ \, x -  v_j \, \Delta t \,$ for $ \, 0 \leq j \leq q-1 . \,$
The time iteration of the scheme can be written as 
\moneq   \label{iteration}
f_j(x, \, t + \Delta t) \,=\,  f_j^* ( x - v_j \, \Delta t ,\, t)  
\,, \qquad 0 \leq j \leq q-1 \,, \quad x \in {\cal L}. 
\monend  

 \monitem   
From the knowledge of the previous  algorithm, it is possible to derive a set of
equivalent partial differential equations for the conserved variables. If the Gaussian 
functions $ \, G_k \,$ are  linear, this set of equations takes the form
\moneq   \label{edp-lineaire}
{{\partial W}\over{\partial t}} \,-\,  \alpha_1 \, W \,-\,  \Delta t \,  \alpha_2 \, W 
\,-\, \cdots \,-\,  \Delta t^{j-1}  \,  \alpha_j \, W \,-\, \cdots \,-\,  \,=\, 0  \,,
\monend  
where $ \, \alpha_j \,$ is for $ \, j \geq 1 \,$ is a space derivation operator of
order $j$. We refer the reader to \cite{Du08}   for the presentation
of our approach in the general case.
In this contribution, we have developed an explicit algebraic linear
version of the algorithm  detailed in Appendix~1. Moreover, 
we consider that the lattice Boltzmann scheme 
is invariant by rotation at order $ \, \ell \,$
if the equivalent partial differential equation 
\moneq   \label{edp-lineaire-ell}
{{\partial W}\over{\partial t}} \,-\, \sum_{j=1}^{\ell} 
 \Delta t^{j-1}  \,  \alpha_j \, W   \,=\, 0   
\monend  
obtained from (\ref{edp-lineaire}) by truncation at the order $ \, \ell \,$
is invariant by rotation. 
For acoustic-type models, the partial differential equation 
(\ref{edp-lineaire-ell}) has to be identical to 
(\ref{rot-eqs-2d}) or (\ref{rot-eqs-3d}) for dimension 2 or~3. 
In the following, we determine the equivalent partial differential equations
for classical lattice Bolzmann schemes in the general linear case. 
Then we fit the equilibrium and relaxation parameters of the scheme             
in order to enforce rotational invariances at all orders between 1 and 4.

  \bigskip \bigskip   \noindent {\bf \large 4) \quad  D2Q9   } 
  
 \noindent
The isotropy of the lattice Boltzmann scheme D2Q9 for the acoustic equations has been 
studied in detail in \cite{ADG12,ADGG13}. We give here only a summary of our results.

\monitem 
The matrix M for the D2Q9 lattice Boltzmann   model is 
of the form
\moneq   \label{d2q9-matrice-M} 
M_{kj} = p_k(v_j) \,, \qquad 0 \leq \, j \,, \,\, k   \, \leq  q-1 \,  
\monend 
with polynomials $ \, p_k \, $  detailed  in the contribution \cite{LD13}.  
With this  choice,   
the  moments are named according to the following notations:   
\moneq   \label{d2q9-moments}  
\left\{ \begin{array}{rcl}
 0 &  1   & \lambda^0  \\
 1,2   &   X, Y    &  \lambda^1 \\
 3 & \varepsilon  &  \lambda^2  \\
 4, 5  &  XX, XY   &  \lambda^2 \\
 6, 7  &  q_x, q_y   &  \lambda^3  \\
 8  & \varepsilon_2   &  \lambda^4 \, . 
   \end{array} \right. \monend  
We have recalled in (\ref{d2q9-moments}) the degrees of 
homogeneity of each moment $ \, m_k \,$ 
relative to the reference numerical velocity
$\, \lambda \,\equiv \, {{\Delta x}\over{\Delta t}} $. 

\monitem 
At  first order, the 
invariance by rotation (\ref{rot-eqs-2d}) takes the form 
\moneq   \label{edp-ordre-1}  
\left\{\begin{array}{l} 
\partial_t \rho +   \,  {\rm div} J  \, = \,  {\rm O}(\Delta t)   \\
\displaystyle \partial_t J + c_0^2 \, \nabla \rho   \, = \,   {\rm O}(\Delta t) 
\end{array}\right. \monend  
if the next moments of degree 2 follow the  simple  equilibrium:     
\moneq   \label{d2q9-cns-ordre-1}  
 \varepsilon^{\rm eq}  = \alpha \, \lambda^2 \, \rho   \,, \quad  
  XX^{\rm eq}  = XY^{\rm eq}    = 0  \, . \monend  
Then the sound velocity $\, c_0 \,$ in the equation (\ref{edp-ordre-1})
satisfies   
\moneq   \label{d2q9-coefs-ordre-1}  
c_0^2 = {{\lambda^2}\over{6}} \, ( 4 + \alpha)   \,.  
\monend

\monitem   
At second  order,  the 
invariance by rotation (\ref{rot-eqs-2d}) is realized under specific  conditions
for the third order moments   $\, q  \equiv (q_x , \, q_y) .\,$ 
This equilibrium condition is defined with the help of 
H\'enon's  \cite{He87}
parameters $\, \sigma_k \,$ defined 
from the coefficients $\, s_k \,$  according to 
\moneq   \label{sigmas}   \sigma_k  
\, \equiv \,    {{1}\over{s_k}} - {{1}\over{2}} \,\qquad  {\rm when} \,\,  k \geq 3 \, . 
 \monend  
The stability condition (\ref{stabi}) can be written as
\moneq   \label{stabi-sigmas} 
\sigma_k \, > \, 0 \,, \quad{\rm for} \,\,  k \geq N \, . 
 \monend  
We have  
\moneq   \label{d2q9-cns-ordre-2}   
 q^{\rm eq}  =
 {{ \sigma_4 - 4 \,  \sigma_5 } \over {  \sigma_4 + 2 \,  \sigma_5 }}
 \, \lambda^2 \, J \,   . 
 \monend  
We obtain with  these conditions   the following  equivalent partial differential
equations 
\moneq   \label{edp-ordre-2}  
\left\{\begin{array}{l}  
\partial_t \rho +   \,  {\rm div} J  \, = \,   {\rm O}(\Delta t^2)   \\
\displaystyle \partial_t J + c_0^2 \, \nabla \rho  - \mu \, \triangle J 
- \zeta \, \nabla \, {\rm div} J \, = \,  {\rm O}(\Delta t^2) \, . 
\end{array}\right. \monend  
The physical 
viscosities $ \, \mu \,$ and  $ \, \zeta \,$ can be determined
according to 
\moneq   \label{d2q9-coefs-ordre-2}  
\mu  =     {{ \sigma_4 \, \sigma_5 }\over { \sigma_4 +  2 \, \sigma_5}}  \, \lambda \,
 \Delta  x   \,, \quad 
   \zeta =  \sigma_3 \, 
{{ ( 2 \,  \sigma_4 -  2 \,  \sigma_5 - \alpha  \,   \sigma_4 - 2 \, \alpha  \,   \sigma_5 ) } 
 \over { 6 \, ( \sigma_4 +  2 \, \sigma_5 ) }} 
 \,  \lambda \, \Delta  x  \, . \monend  
We observe that the classical isotropy condition $\,  \sigma_4 =  \sigma_5 \,$ for
the second order moments $ \, XX \,$ and $ \, XY \,$ is not necessary for  the 
relaxation at this  second order  step. 

\monitem  
At third order,  the system  (\ref{rot-eqs-2d}) takes the particular form
\moneq   \label{edp-ordre-3}  
\left\{\begin{array}{l}  
\partial_t \rho +   \,  {\rm div} J  + \xi \, \triangle \,  {\rm div} J  \, = \,  
  {\rm O}(\Delta t^3)   \\
\displaystyle \partial_t J + c_0^2 \, \nabla \rho  - \mu \, \triangle J 
- \zeta \, \nabla \, {\rm div} J + \chi \, \nabla \triangle  \rho  
\, = \,  {\rm O}(\Delta t^3) \, . 
\end{array}\right. \monend  
This is possible if the complementary relations
\moneq   \label{d2q9-cn-ordre-3}   
\sigma_4 \, = \,  \sigma_5  \,,  \quad 
\varepsilon_2 ^{\rm eq}  = - {{\lambda^4 \, \rho}\over{2}} \, \big( 
3 \, \alpha +4 \big)  
\monend  
hold. Then the heat flux at equibrium has an expression  (\ref{d2q9-cns-ordre-2})
that can be simply written as
\moneq   \label{d2q9-cn2-ordre-3}
 q^{\rm eq}  = - \, \lambda^2 \, J   \, . 
 \monend  
The coefficients in the equations 
(\ref{edp-ordre-3}) are given by the expressions 
\moneq   \label{d2q9-coefs-ordre-3}   
\left\{\begin{array}{l}  \displaystyle  
\mu  = {{ 1  }\over { 3 }}  \, \sigma_4 \,  \lambda \, \Delta  x  \, , \quad  
\zeta = - {{1}\over{6}} \, \sigma_3 \, \alpha  \,  \lambda \, \Delta  x  \,, \quad 
\xi  \,=\, {{1}\over{72}} \, ( \alpha - 2 ) \, \Delta x^2 \,,    \\ 
 \displaystyle \qquad   \chi  = 
 {{1}\over{216}} \, ( \alpha + 4 ) \, \big(  2 + 6 \, \alpha \, \sigma_3^2 - \alpha 
- 12 \, \sigma_4^2 \big)  \, \lambda^2 \,  \Delta x^2  \, . 
\end{array}\right. \monend   
We remark that the dissipation of the acoustic waves 
$ \, \gamma \equiv {{\mu + \zeta}\over{2}} \,$ (see {\it e.g.}  Landau and Lifshitz
\cite{LL59})  is given according to 
\moneq   \label{d2q9-dissip-son}   
\gamma \,=\,  {{\lambda \,     \Delta x}\over{12}} \, 
\big( 2 \sigma_4 - \alpha \, \sigma_3 \big) \,     . 
\monend   

\monitem  
The invariance  by rotation  at fourth  order of the D2Q9 lattice 
Boltzmann scheme does not give completely  satisfactory results, as observed  
previously in  \cite{ADGG13}.                                                  
In order to get equivalent equations at order~4 of the type   
\moneq   \label{edp-ordre-4}  
\left\{\begin{array}{l}  
\partial_t \rho +   \,  {\rm div} J  + \xi \, \triangle \,  {\rm div} J 
+ \eta \, \triangle^2 \rho \, = \,   {\rm O}(\Delta t^4)   \\
\displaystyle \partial_t J + c_0^2 \, \nabla \rho  - \mu \, \triangle J 
- \zeta \, \nabla \, {\rm div} J + \chi \, \nabla \triangle  \rho  
+  \, \mu_4  \, \triangle^2  J 
  +  \zeta_4  \, \nabla \, {\rm div} \, \triangle \, J  =   {\rm O}(\Delta t^4)   ,  
\end{array}\right. \monend  
%
%
it is necessary to fix some relaxation parameters:
\moneq   \label{d2q9-cns-ordre-4}  
\sigma_3 = \sigma_4 =  \sigma_8   \,, \quad 
 \sigma_6 = \sigma_7 = {{1}\over{6 \, \sigma_4}}  \, . 
 \monend  
The two viscosities $ \mu   $ and $ \zeta $ are now   dependent and 
the dissipation of the acoustic waves introduced previously in 
(\ref{d2q9-dissip-son}) admits now the expressions 
\moneq   \label{d2q9-mu-zeta-ordre-4}   
\mu  = {{ 1  }\over { 3 }}  \, \sigma_4 \,  \lambda \, \Delta  x \,, \quad  
\zeta = - {{1}\over{6}} \, \sigma_4 \, \alpha  \,  \lambda \, \Delta  x  \,, \quad  
\gamma \,=\,  {{\lambda \,  \sigma_4 \,     \Delta x}\over{12}} \, (2 - \alpha) \,. 
 \monend  
Observe that the dissipation $\, \gamma \,$  
 is positive under usual conditions. 
If the conditions (\ref{d2q9-cns-ordre-4}) are satisfied, we can specify   
the   coefficients of the fourth order terms in the equations (\ref{edp-ordre-4}): 
\moneq   \label{d2q9-coefs-ordre-4}   
\left\{\begin{array}{l}  \displaystyle
 \eta \,= \, {{\lambda \,  \Delta x^3}\over{432}} \, 
( \alpha + 4) \, ( \alpha - 2 ) \, , \quad  
\mu_4  \,=\,  {{\lambda \,  \Delta x^3}\over{108}} \, 
(12 \, \sigma_4^2 -1 )  \, , \quad       
  \\  \vspace  {-.4 cm}   \\      \displaystyle    
  \zeta_4  \,=\,   {{\lambda \,  \Delta x^3}\over{216}} \, 
\sigma_4 \, \Big( 12 - \alpha - 2 \, \alpha^2 
+ 12 \, \sigma_4^2 \, ( \alpha^2  - \alpha -4 ) \Big)   \, . 
\end{array}\right. \monend   
%
%

\monitem  
In \cite{ADGG13}, we have conducted  a set of numerical experiments that 
  make more explicit    
the isotropy qualities of four different variants of the D2Q9 lattice
Boltzmann scheme for the numerical simulation of 
acoustic waves. The   orders of isotropy precision   numerically computed
are in coherence with the  level of accuracy  presented in this section.

  \bigskip \bigskip   \noindent {\bf \large 5) \quad  D2Q13   } 

\monitem
Four more velocities are added to the D2Q9 scheme to construct D2Q13.
The details can be found {\it e.g.} in   \cite{LD13}.                         
Nine  moments are analogous to those proposed  in (\ref{d2q9-moments})        
for  D2Q9 and four moments ($r_x$, $r_y$, $\varepsilon_3$, $XX_e$) are new: 
%
\moneq   \label{d2q13-moments}  
\left\{ \begin{array}{rcl}  
 0 &  1   & \lambda^0  \\
 1,2   &   X, Y    &  \lambda^1 \\
 3 & \varepsilon  &  \lambda^2  \\
 4, 5  &  XX, XY   &  \lambda^2 \\
 6, 7  &  q_x, q_y   &  \lambda^3  \\
 8, 9  &  r_x, r_y   &  \lambda^5  \\
 10 & \varepsilon_2   &  \lambda^4 \\
 11 & \varepsilon_3   &  \lambda^6 \\
 12  &  XX_e    &  \lambda^4 \, . 
   \end{array} \right. \monend  
%

\monitem 
At  first order, the 
invariance by rotation (\ref{rot-eqs-2d}) takes again the form 
 (\ref{edp-ordre-1}). 
The conditions (\ref{d2q9-cns-ordre-1}) are essentially unchanged, 
except that the sound velocity is now evaluated according to the relation  
\moneq   \label{d2q13-cns-ordre-1}   
c_0^2 \,=\,  {{\lambda^2}\over{26}} \, ( 28 + \alpha) \,\equiv\, c_s^2 \,\lambda^2   \,.
  \monend
%
\monitem 
At second  order, the equivalent partial differential equations 
take the  isotropic form  (\ref{edp-ordre-2})   
when we  have : 
\moneq   \label{d2q13-cns-ordre-2}   
q ^{\rm eq}  = \varphi \, \lambda^2 \, J \,, \quad   
r ^{\rm eq}  = {{1}\over{12}} \,  \Big( {{ 20 \,\sigma_5  
- 85 \, \sigma_4 - 49 \, \varphi \,  \sigma_4 
- 14 \, \varphi \,  \sigma_5   }\over{ \sigma_4 + \sigma_5}}   \Big) 
 \, \lambda^4 \, J \, . 
 \monend  
Then the isotropy coefficients in  (\ref{edp-ordre-1}) 
have the following expressions:
\moneq   \label{d2q13-coefs-ordre-2}  
\mu  =   {{\lambda \, \Delta  x}\over{2}} \,   {{ \sigma_4 \, \sigma_5 }\over { \sigma_4 +  \sigma_5}}  \,
( 3 + \varphi) \,     \,, \quad 
  \zeta =  {{\lambda \, \Delta  x}\over{26}} \, \sigma_3 \, 
( 11 + 13 \, \varphi - \alpha )   \,, \quad  
 \sigma_k > 0 \, \, \, {\rm when} \,\,  k \geq 3 \, .  
 \monend

\monitem 
The invariance by rotation  at third order  of the mass equation is realized if
we impose a unique value for the  relaxation coefficients of the second
order moments $XX$ and $XY$ introduced in  (\ref{d2q13-moments}): 
\moneq   \label{d2q13-cn-1-ordre-3} 
\sigma_4 \,=\, \sigma_5 \, . 
\monend  
Moreover, the attenuation of sound  waves $ \, \gamma \,$ does   not 
 depend at first order   on the advective velocity if 
(\ref{d2q13-cn-1-ordre-3})  is satisfied. 
For the 
invariance by rotation  of the momentum equation, we must impose also 
the following equilibrium values 
 for the moments $\, m_{10} \equiv \varepsilon_2 ,\, $  
$ m_{11} \equiv \varepsilon_3 \, $ and  $\, m_{12} \equiv  XX_e $:
\moneq   \label{d2q13-cn-2-ordre-3} 
\left\{\begin{array}{l}  \displaystyle  
\varepsilon_2 ^{\rm eq}  = \Big(  -5\, \alpha  
+ {{77}\over{26}} \, \, \varphi \, \alpha \, 
 +  {{1078}\over{13}}  \, \varphi \Big)  \, \lambda^4 \, \rho   
  \\  \vspace  {-.4 cm}   \\      \displaystyle   
\varepsilon_3 ^{\rm eq}  =  \Big(  {{\alpha}\over{ 48}}  - {{137 }\over{ 12}}  
- {{135 }\over{ 208}} \, \alpha \, \varphi  - {{945 }\over{ 52}} \, \varphi  \Big)  \, 
 \lambda^4 \, \rho   \, , \qquad    XX_e^{\rm eq} \,=\, 0 \, . 
\end{array}\right. \monend  
Then  the equivalent equations at order~3 of the D2Q13 lattice Boltzmann scheme
are still given by the equations (\ref{edp-ordre-3}). 
The equilibrium condition (\ref{d2q13-cns-ordre-2}) is  now  written as 
\moneq   \label{d2q13-cn-3-ordre-3} 
r ^{\rm eq}  =  - {{1}\over{24}} \, ( 65 + 63 \,  \varphi )  
 \, \lambda^4 \, J \,   
 \monend 
and  the coefficients in the equations  (\ref{edp-ordre-3}) 
can be clarified:     
\moneq   \label{d2q13-coefs-ordre-3} 
\left\{\begin{array}{l}  \displaystyle    
\mu  = {{ 1  }\over { 4 }}  \, \sigma_5 \, ( 3 + \varphi) \, \lambda \, \Delta  x \,, \quad 
 \zeta =  {{1}\over{26}} \, \sigma_4 \, 
( 11 + 13 \, \varphi - \alpha ) \,  \lambda \, \Delta  x   \,,
  \\  \vspace  {-.5 cm}   \\      \displaystyle     
 \xi  \,=\, {{1}\over{624}} \, ( 2 \, \alpha - 39 \, \varphi - 61) \, \Delta x^2  \,, 
  \\  \vspace  {-.5 cm}   \\      \displaystyle    
 \chi \, = \,  {{1}\over{8112}}  \, (28 + \alpha ) \, 
  \Big( 61  +   39  \, \varphi  + 12 \, \alpha \, \sigma_3^2
- 2 \, \alpha - 78 \, \varphi \, \sigma_4^2   
  \\  \vspace  {-.5 cm}   \\      \displaystyle  
  \qquad  \qquad  \qquad  \qquad  \qquad  \qquad  \qquad 
- 156  \, \varphi \, \sigma_3^2 
- 234  \, \sigma_4^2 - 132 \,  \sigma_3^2
\Big)  \, \, \lambda^2 \,  \Delta x^2 \, . 
\end{array}\right. \monend

\monitem  
The invariance   by rotation  at fourth  order 
is satisfied if we add to the previous conditions 
(\ref{d2q13-cn-1-ordre-3})  (\ref{d2q13-cn-2-ordre-3}) and  
(\ref{d2q13-cn-3-ordre-3}) the new ones:
\moneq   \label{d2q13-cn-1-ordre-4}  
\left\{\begin{array}{l}  \displaystyle 
q ^{\rm eq}  \,=\,  - {{7}\over{5}}  \, \lambda^2 \, J  \,, \quad 
  \sigma_6 = \sigma_7  =  {{1}\over{12 \, \sigma_4}}  \,,  
  \\  \vspace  {-.5 cm}   \\      \displaystyle  
\sigma_8   \,=\, \sigma_9   \,=\, {{5}\over{24}} \,\, {{155-a}\over{a -308}} \,\, {{1}\over{\sigma_4}} 
\,+\,  {{1}\over{24}} \,\, {{7\, a - 1391}\over{a-308}} \,\, {{1}\over{\sigma_3}}  \,,
  \\  \vspace  {-.5 cm}   \\      \displaystyle  
\sigma_{10} \,=\,  {{3973}\over{45}} \,\, {{43\, a - 16610}\over{89 \, a- 20680}} \, \,
 {{5 \, c_s^2-4}\over{1189\, c_s^2 - 828}} \,\, \sigma_3 \,\,+\, 
  \\  \vspace  {-.5 cm}   \\      \displaystyle    \qquad  \qquad  \qquad  \qquad  \qquad 
\,+\,  {{154}\over{1395}} \,\, {{7\, a - 1391}\over{89 \, a - 20680}} \, \,
 {{725 \, c_s^2-418}\over{1189\, c_s^2 - 828}} \,\, a \, \sigma_4  \,,
  \\ \displaystyle    
\sigma_{11} \,=\,  {{a}\over{155}} \,\, \sigma_4 \, . 
\end{array}\right. \monend   
If the parameters $\, c_s \,$ and  $\, \alpha \,$ relative to the non-dimensionalized sound
velocity are  linked together thanks to (\ref{d2q13-cns-ordre-1})    
and if the new parameter $\, a \,$ is  chosen such that 
\moneq   \label{d2q13-cn-3-ordre-4}  
c_s^2   \,< \, {{418}\over{25}} \,, \quad 
-28   \,< \, \alpha \leq  - {{9432}\over{725}} \simeq -13 \,, \quad 
155 \,< \, a  \,< \, {{1391}\over{7}}  \simeq 198 \,, 
 \monend 
the coefficients $\, \sigma_8, $  $\, \sigma_{10} \, $ and  
$\, \sigma_{11} \, $ are strictly positive if it is the case for 
$\, \sigma_{3} \, $ and $\, \sigma_{4} . \, $
In this case, the stability conditions (\ref{stabi-sigmas})
are satisfied for the coefficients $\, \sigma_3,$ $\, \sigma_4,$ $\, \sigma_8,$
$\, \sigma_{10} $ and $\, \sigma_{11}$.    
%
With the choice (\ref{d2q13-cn-1-ordre-4})  
the nontrivial algebraic expressions of the previous  conditions 
(\ref{d2q13-cn-1-ordre-3}), (\ref{d2q13-cn-2-ordre-3}) and  
(\ref{d2q13-cn-3-ordre-3}) can be written as 
\moneq   \label{d2q13-cn-4-ordre-4}  
  r ^{\rm eq}  =  {{29}\over{30}}  \, \lambda^4 \, J \,, \quad 
 \varepsilon_2 ^{\rm eq}  =  - \Big(
 {{1189}\over{130}}\, \alpha + {{7546}\over{65}} \Big) \,  \lambda^4 \, \rho \,, \quad  
\varepsilon_3 ^{\rm eq}  = 
 \Big( {{547}\over{39}} + {{145}\over{156}}\, \alpha \Big) \, \lambda^4 \, \rho \, . 
 \monend 
With the above conditions  
(\ref{d2q13-cn-1-ordre-4})   and  (\ref{d2q13-cn-4-ordre-4}) 
the equivalent equations of the D2Q13 lattice  Boltzmann scheme at fourth order 
are  made explicit     
    in  (\ref{edp-ordre-4}), with the associated coefficients, except
$ \, \zeta_4$, given according to: 
\moneq   \label{d2q13-coefs-ordre-4} 
\left\{\begin{array}{l}  \displaystyle   
\xi \,=\, {{5\, \alpha-16}\over{1560}}  \, \Delta x^2 \,, \quad 
\eta \,=\,  {{\alpha + 28}\over{40560}}   \, \big( 
36 \sigma_3 + 5 \sigma_3 \, \alpha - 52 \sigma_4 \big) \, \lambda \, \Delta x^3  
  \\  \vspace  {-.5 cm}   \\      \displaystyle   
\mu  = {{ 2  }\over { 5 }}  \, \sigma_4 \, \lambda \, \Delta  x   \,, \quad 
 \zeta =     - {{1}\over{130}} \, \sigma_3 \, 
( 36 + 5 \,  \alpha ) \,  \lambda \, \Delta  x    
  \\  \vspace  {-.5 cm}   \\      \displaystyle  
\chi  \,=\,  {{1}\over{20280}} \, (28 + \alpha) \, 
  \Big( 16  - 5 \, \alpha     + 216 \, \sigma_3^2  - 312  \, \sigma_4^2 
+ 30 \, \alpha \, \sigma_3^2  \Big)   \,\, \lambda^2 \,  \Delta x^2 
  \\  \vspace  {-.5 cm}   \\      \displaystyle   
\mu_4  = {{ \sigma_4 \, \lambda \, \Delta  x^3 }\over { 300\, \sigma_3 \, (a - 308)  }}  \, 
\Big(  4483 \sigma_4 - 5099 \, \sigma_3 - 23 \, a \, \sigma_4
+ 25 \, a \, \sigma_3 
\\      \displaystyle  \qquad  \qquad   \qquad    \qquad  \qquad   \qquad  
- 14784 \, \sigma_3 \, \sigma_4^2  + 48 \, a \,  \sigma_3 \, \sigma_4^2 \Big)   \, . 
\end{array}\right. \monend 
The algebraic expression of the coefficient $ \, \zeta_4 \,$ 
is quite long.  
With   H\'enon's coefficients $ \, \sigma_j \, $ defined according to (\ref{sigmas}), 
the related moments numbered by the relations (\ref{d2q13-moments}),
the equilibrium  of the energy   (\ref{d2q9-cns-ordre-1}) parametrized by $\, \alpha $, 
and the parameter $a$ introduced at (\ref{d2q13-cn-1-ordre-4}), 
the coefficient $\, \zeta_4 \, $
for the fourth order term in (\ref{edp-ordre-4}) can be evaluated according to:

\moneqstar 
\left\{     \begin{array}{l}  \displaystyle  
\zeta_4  = {{ \sigma_4 \, \lambda \, \Delta  x^3 }\over { 56581200 \,\, \sigma_3 \,
\, ( 89 a - 20680)  }}  \, \Big( 525433428 \, a \, \sigma_3 \, \sigma_4  
\,+\,  576972000 \, \, \alpha^2  \sigma_3^2  
\\      \displaystyle  \qquad  
+\, 18001526400 \,  \alpha  \, \sigma_3^3 \, \sigma_4  
\,+\, 18001526400 \,  \alpha \,   \sigma_3 \,   \sigma_4^3 
+\, 65975 \, a^2 \,   \alpha \, \sigma_3   \, \sigma_4
\\      \displaystyle  \qquad  
+\, 170558856 \, a \, \sigma_4^2 
 \,-\, 334055628 \, a \, \sigma_3^2 
\,-\, 858312 \, a^2 \, \sigma_4^2 
\,-\, 159243217380 \, \sigma_3 \, \sigma_4
\\      \displaystyle  \qquad   
+\, 75143778660 \, \sigma_3^2 \, 
+\, 18001526400 \, \alpha  \, \sigma_3^2   \, \sigma_4^2  
\,-\, 77472720 \,   a \, \alpha  \, \sigma_3^3  \, \sigma_4  
\\      \displaystyle    \qquad 
-\, 77472720 \, a \, \alpha  \,   \sigma_3^2   \,  \sigma_4^2  
\,-\, 77472720  \, a \, \alpha   \,  \sigma_3 \,     \sigma_4^3 
+\, 504042739200  \, \sigma_3  \, \sigma_4^3 
\\      \displaystyle    \qquad 
+\, 129610990080 \, \sigma_3^3 \, \sigma_4  
\,+\, 129610990080 \, \sigma_3^2  \, \sigma_4^2  
+\, 858312 \, a^2  \, \sigma_3 \, \sigma_4 
\\      \displaystyle    \qquad 
+\, 22940190 \, a \, \alpha  \,   \sigma_3  \, \sigma_4   
\,+\, 17841109925  \, \alpha  \, \sigma_3^2 
\,-\, 65975  \, a^2 \, \alpha \,  \sigma_4^2  
\\      \displaystyle    \qquad 
+ \, 13110175 \, a  \, \alpha  \, \sigma_4^2   
\,-\, 3461832000 \, \alpha^2  \, \sigma_3^4  
\,-\, 85853433600 \, \alpha  \, \sigma_3^4   
\\      \displaystyle    \qquad 
-\, 2483100 \, a  \, \alpha^2   \, \sigma_3^2  
\,-\, 78263065 \,  a \, \alpha    \, \sigma_3^2  
\,-\, 438683351040 \, \sigma_3^4
\\      \displaystyle    \qquad  
+\, 369485280 \, a \, \alpha  \, \sigma_3^4  
\,+\, 14898600 \, a  \, \alpha^2 \, \sigma_3^4    
\, +\, 1887950592 \, a  \, \sigma_3^4  
\\      \displaystyle    \qquad  
-\, 557803584 \, a \, \sigma_3^2 \, \sigma_4^2 
\,-\, 2169236160 \, a \, \sigma_3   \, \sigma_4^3  
\,-\, 557803584  \, a \, \sigma_3^3 \, \sigma_4   
\\      \displaystyle    \qquad  
-\, 8032585925 \,  \alpha  \,  \sigma_3 \, \sigma_4  \Big) \, . 
\end{array}\right. \monendstar

  \bigskip \bigskip  \newpage
 \noindent {\bf \large 6) \quad  D3Q19   } 

\monitem 
For the scheme  D3Q19, 
we have  4 conservation laws and a total of 19 moments. 
We refer  {\it e.g.} to  \cite{DL09}  
for an algebraic expression of the polynomials $ \, p_k \, $ in (\ref{d2q9-matrice-M}):
\moneq   \label{d3q19-moments}  
\left\{ \begin{array}{rcl}  
 0 &  1   & \lambda^0  \\
 1, 2, 3   &   X, Y, Z    &  \lambda^1 \\
 4 & \varepsilon  &  \lambda^2  \\
 5, 6   &  XX, WW   &  \lambda^2 \\
 7, 8, 9  &  XY, YZ, ZX   &  \lambda^2 \\
 10, 11, 12  &  q_x, q_y , q_z   &  \lambda^3  \\
 13 & \varepsilon_2   &  \lambda^4 \\
 14, 15  &  XX_e, WW_e   &  \lambda^4 \\
 16, 17, 18  &  {\rm antisymmetric \ of \ order} \  3   &  \lambda^3 \, . 
   \end{array} \right. \monend  
The results summarized in this Section have been essentially      
considered (quickly) in the previous contribution 
\cite{DL11}. They have been also used by 
Leriche, Lallemand and Labrosse \cite{LLL08} for 
the numerical determination of the eigenmodes of the Stokes problem 
in a cubic cavity.

\monitem 
We write the four equivalent partial differential equations at first order with the method
explained in Appendix~1. 
Then we impose that the associated modes are isotropic, {\it i.e.} contain only
partial differential operators that are invariant by rotation. 
In Fourier space, the coefficients of the associated determinant must 
contain only powers of the wave vector. 
The associated equations are highly nonlinear relative to the 
coefficients of the equilibrium matrix introduced in (\ref{gaussian}). 
We have obtained a family of parameters  by enforcing 
the linearity of the solution of the isotropic equations. 
With this  constraint, we have to  impose a relation for the ``energy'' moment 
 $\, m_4 \equiv \varepsilon \,$ at equilibrium:  
\moneq   \label{d3q19-cn-1-ordre-1}   
 \varepsilon^{\rm eq}  \,=\,  \alpha \, \lambda^2 \, \rho   \, . 
\monend  
Moreover, the equilibrium values
for the  moments $\, m_5 \,$ to  $\, m_9 \,$   of degree two introduced in 
(\ref{d3q19-moments}) are equal to zero: 
\moneq   \label{d3q19-cn-2-ordre-1}   
XX^{\rm eq}  \,=\,   WW^{\rm eq}    \,=\,  XY^{\rm eq}     
\,=\,  YZ^{\rm eq}    \,=\,   ZX^{\rm eq}   \,=\,  0  \, . 
\monend  
When the conditions (\ref{d3q19-cn-1-ordre-1}) and  (\ref{d3q19-cn-2-ordre-1})
are realized, the isotropic equivalent system is given 
by the system of first order acoustic equations  (\ref{edp-ordre-1}). 
Moreover, the sound velocity $ \, c_0 \,$ satisfies 
\moneq   \label{d3q19-coefs-ordre-1} 
 c_0^2 \,=\, {{\alpha + 30}\over{57}} \, \lambda^2 \,\, \equiv \,\, c_s \,  \lambda^2 \, . 
\monend  
%

\monitem 
Invariance by rotation  at second  order is realized if we impose on one hand  
\moneq   \label{d3q19-cns-1-ordre-2}   
q ^{\rm eq}  = 2 \, {{3 \, \sigma_5 - 4 \, \sigma_7 }\over{ \sigma_5 + 2 \  \sigma_7 }}   \, 
\lambda^2 \, J \,  
 \monend  
and on the other hand 
\moneq   \label{d3q19-cns-2-ordre-2}   
\sigma_5 \, = \, \sigma_6 \,, \qquad \sigma_7 \,=\,  \sigma_8 \, = \,  \sigma_9 \, . 
 \monend  
Then the equivalent partial differential equations 
of the D3Q19 lattice Boltzmann scheme take the form  (\ref{edp-ordre-2}). 
The associated coefficients are given according to 
\moneq   \label{d3q19-coefs-ordre-2} 
\left\{\begin{array}{l}  \displaystyle   
\mu  = {{ \sigma_5 \, \sigma_7 }\over { \sigma_5 + 2 \, \sigma_7}}  \, \lambda \, \Delta  x 
  \\  \vspace  {-.5 cm}   \\      \displaystyle   
 \zeta =  {{\lambda \, \Delta  x}\over{  57 \, (  \sigma_5 + 2 \, \sigma_7 ) }}
\, \big(  27 \, \sigma_4 \,  \sigma_5 \,+\, 19 \, \sigma_5 \, \sigma_7  
\,  - \, 22  \, \sigma_4 \, \sigma_7  \,-\,  \alpha \,\sigma_4 \, \sigma_5 \, 
\,- \, 2 \, \alpha   \, \sigma_7 \, \sigma_4 \, \alpha   \big)  \,  \, . 
   \end{array} \right. \monend  
%

\monitem 
At third order, 
if we impose the previous relations   
(\ref{d3q19-cn-1-ordre-1})  (\ref{d3q19-cn-2-ordre-1})  
 (\ref{d3q19-cns-1-ordre-2}) and     (\ref{d3q19-cns-2-ordre-2}),  {\it id est} 
an equilibrium for the ``energy square'' $ \, \varepsilon_2 \,$
given below,  a null value for $ \, m_{14}^{\rm eq}   \,$ and  $ \, m_{15}^{\rm eq}   \,$ 
and a supplementary condition for the relaxation coefficients, 
{\it id est} 
\moneq   \label{d3q19-cn-1-ordre-3}   
 \varepsilon_2^{\rm eq}  \,=\,  {{42 + 9 \, \alpha}\over{9}} \, \lambda^2 \, \rho  \, , \qquad 
 XX_e^{\rm eq} \, = \,  WW_e^{\rm eq}  \,=\, 0     \,, \qquad 
\sigma_5 \,=\, \sigma_7 \,, 
\monend  
the equivalent equations of the D3Q19 lattice Boltzmann scheme 
are exactly given by   (\ref{edp-ordre-3}). 
We observe that the relation  (\ref{d3q19-cns-1-ordre-2})   takes now the form
\moneq   \label{d3q19-cn-2-ordre-3}   
q ^{\rm eq}  = -  \, {{2}\over{3}}   \, \lambda^2 \, J  \,  
\monend  
and the coefficients  associated to the equations  (\ref{edp-ordre-3}) 
 can be deduced through  an elementary process:
\moneq   \label{d3q19-coeffs-ordre-3}  
\left\{\begin{array}{l}  \displaystyle  
\mu  = {{ 1  }\over { 3 }}  \, \sigma_5 \,  \lambda \, \Delta  x  \, , \quad  
\zeta = {{\lambda \, \Delta  x }\over{171}}   
\, \big(  5 \, \sigma_4  \, + \, 19 \, \sigma_5 \,-\, 3 \, \alpha \, \sigma_4   \big) \,, \quad   
\xi  \,=\, {{\Delta x^2 }\over{684}} \, ( \alpha - 27 )  \,, 
  \\  \vspace  {-.5 cm}   \\      \displaystyle   
  \chi  \,=\,  {{ \lambda^2 \,  \Delta x^2 }\over{19494}} \,  \, \big( \alpha+30 \big) \, 
\big( 27 \,+\, 6 \, \alpha \, \sigma_4^2   \,-\, 10\, \sigma_4^2   
\,-\, 152 \, \sigma_5^2  \, -\, \alpha  \big)  \, . 
\end{array}\right. \monend   
%
%
The bulk viscosity $ \, \zeta_b \equiv 3 \, \zeta - \mu \, $   
(see {\it  e.g.}  Landau and Lifshitz \cite{LL59})  is essentially function of  
the relaxation parameter associated to the energy $ \, \varepsilon \,$: 
\moneq   \label{d3q19-bulk-ordre-3}  
\zeta_b  \,= \,  {{\lambda \, \sigma_4  \,  \Delta  x }\over{57}} \, \big( 5 - 3 \, \alpha \big) \, .   
\monend   
%

\monitem 
We have found also a variant of the previous relations to enforce third order    
isotropy. We can replace the relations (\ref{d3q19-cn-1-ordre-3}) by the 
following ones, with an  undefined parameter $ \, \beta $:  
\moneq   \label{d3q19-cn-3-ordre-3}   
\left\{\begin{array}{l}  \displaystyle
 \varepsilon_2^{\rm eq}  \,=\,  \beta \, \lambda^2 \, \rho  \, , \quad 
 XX_e^{\rm eq} \, = \,  WW_e^{\rm eq}  \,=\, 0     \,,  \quad 
\sigma_5 \,=\, \sigma_7  \,,
  \\  \vspace  {-.5 cm}   \\      \displaystyle   
\sigma_{10} \,=\, \sigma_{11} \,=\, \sigma_{12} \,=\, {{1}\over{12 \, \sigma_5}} \,, \quad 
\sigma_{16} \,=\,  {{1}\over{8 \, \sigma_5}}  \,, \quad 
\sigma_{17} \,=\,  {{1}\over{4 \, \sigma_5}}  \,, \quad 
\sigma_{18} \,=\,  {{1}\over{12 \, \sigma_5}}  \,. 
   \end{array} \right. \monend  
The relations  (\ref{d3q19-cn-2-ordre-3}),     (\ref{d3q19-coeffs-ordre-3}) 
 and (\ref{d3q19-bulk-ordre-3})  are 
not changed, except  that the coefficient  $ \, \chi \, $ in the second line 
of  (\ref{edp-ordre-3}) is now given by
\moneq   \label{d3q19-chi-3-ordre-3}   
\left\{\begin{array}{l}  \displaystyle 
\chi \, = \, {{\lambda^2 \, \Delta x^2}\over{409374 \,\, \sigma_5}} 
\, \Big(
16212 \, \sigma_5 \,+\, 126 \, \alpha^2 \, \sigma_4^2 \, \sigma_5  
\,+\, 3570 \, \alpha \,  \sigma_4^2 \, \sigma_5 
\, - \, 6300 \, \sigma_4^2 \, \sigma_5 
  \\  \vspace  {-.5 cm}   \\      \displaystyle   \qquad  \qquad  \qquad   \qquad
\, - \, \sigma_5 \, \alpha^2
\,-\, 234 \,  \alpha  \, \sigma_5 \,
\,+\, 361 \, \beta \, \sigma_4 \,+\, 171 \, \alpha \, \sigma_4 
  \\  \vspace  {-.5 cm}   \\      \displaystyle     \qquad  \qquad  \qquad   \qquad
\,+\, 798 \, \sigma_4
\, - \, 3192 \, \alpha \, \sigma_5^3  \, - \, 95760 \, \sigma_5^3
\, -\, 361 \, \beta \, \sigma_5  \Big) \, . 
   \end{array} \right. \monend  
%

\monitem 
The invariance by rotation  at fourth  order has   also been considered.
But due to the low number of remaining parameters, 
the family of Boltzmann schemes that we have obtained impose
constraints between physical parameters. 
We must have in particular 
\moneq   \label{d3q19-cns-1-ordre-4}   
\sigma_4 \,=\, \sigma_5 
\monend  
and this relation induces some {\it a priori}  relationship between the     
shear viscosity $ \, \mu \,$  of relation (\ref{d3q19-coefs-ordre-2}) 
and the bulk viscosity  $ \, \zeta_b  \,$ presented
in (\ref{d3q19-bulk-ordre-3}). Moreover, 
the relations  (\ref{d3q19-cn-1-ordre-3})  (\ref{d3q19-cn-2-ordre-3}) must be satisfied and 
  H\'enon's parameters of the multiple relaxation times have to follow, 
with an ordering proposed in  (\ref{d3q19-moments}), 
the  complementary conditions 
\moneq   \label{d3q19-cns-2-ordre-4}  
\sigma_{10} \,=\, \sigma_{11} \,=\, \sigma_{12} \,=\,  {{1}\over{6 \, \sigma_5}}  \,, \quad 
\sigma_{13} \,=\, \sigma_{14} \,=\, \sigma_{15} \,=\, \sigma_5  \,, \quad 
\sigma_{16} \,=\, \sigma_{17} \,=\, \sigma_{18} \,=\,  {{1}\over{6 \, \sigma_5}} \, . 
\monend  
Then the fourth order isotropic equivalent equations (\ref{edp-ordre-4}) are
satisfied and the associated coefficients can be clarified:  
\moneq   \label{d3q19-coefs-ordre-4} 
\left\{\begin{array}{l}  \displaystyle   
\xi \,=\, {{\alpha-27}\over{684}}  \, \Delta x^2 \,, \,\,  
\eta \,=\,  {{ (\alpha + 30) \,  (\alpha - 27) }\over{38988  }}  
\, \sigma_5 \, \lambda \, \Delta x^3  \,, \,\,  
\mu_4 \,=\,  {{ 12 \, \sigma_5^2 - 1 }\over{108}}  \, \sigma_5 \, \lambda \, \Delta x^3  
  \\  \vspace  {-.5 cm}   \\      \displaystyle  
\zeta_4 \,=\,   {{\sigma_5 \, \lambda \, \Delta x^3}\over{38988  }}  \, 
\Big( \, 2062 + 45 \, \alpha - 4 \, \alpha^2 - 5304 \, \sigma_5^2 
- 612 \, \alpha \, \sigma_5^2 + 24 \, \alpha^2 \, \sigma_5^2  \, \Big) \, . 
   \end{array} \right. \monend

  \bigskip \bigskip   \noindent {\bf \large 7) \quad  D3Q27   } 

\monitem 
For the schemes D3Q27, we refer  {\it e.g.} to    \cite{DL11}  
for an algebraic expression of the polynomials $ \, p_k \, $ 
of the relation  (\ref{d2q9-matrice-M}). 
The moments follow now the nomenclature 
\moneq   \label{d3q27-moments}  
\left\{ \begin{array}{rcl}  
 0 &  1   & \lambda^0  \\
 1, 2, 3   &   X, Y, Z    &  \lambda^1 \\
 4 & \varepsilon  &  \lambda^2  \\
 5, 6   &  XX, WW   &  \lambda^2 \\
 7, 8, 9  &  XY, YZ, ZX   &  \lambda^2 \\
 10, 11, 12  &  q_x, q_y , q_z   &  \lambda^3  \\
 13, 14, 15  &  r_x, r_y , r_z   &  \lambda^5  \\
 16 & \varepsilon_2   &  \lambda^4 \\
 17 & \varepsilon_3   &  \lambda^6 \\
 18, 19  &  XX_e, WW_e   &  \lambda^4 \\
 20, 21, 22   &   XY_e, YZ_e, ZX_e   &  \lambda^4 \\
 23, 24, 25  &  {\rm antisymmetric \ of \ order} \  3   &  \lambda^3 \\
 26 &  XYZ   &  \lambda^3  \, . 
   \end{array} \right. \monend

\monitem 
At first order, we follow the same methodology as       
the one presented  for the previous schemes. 
We keep the relation  (\ref{d3q19-cn-1-ordre-1}) for the momentum
 $\, m_4 \equiv \varepsilon \,$ at equilibrium.  
As for the D3Q19 scheme, the  equilibrium values
for the  moments $\, m_5 \,$ to  $\, m_9 \,$   of degree two introduced in 
(\ref{d3q27-moments}) are equal to zero and the relation 
(\ref{d3q19-cn-2-ordre-1}) still holds. Then 
 the first order isotropic equivalent system is given by  (\ref{edp-ordre-1}).
Observe that  the sound velocity $ \, c_0 \,$ now satisfies 
\moneq   \label{d3q27-coefs-ordre-1} 
 c_0^2 \,=\, {{\alpha + 2}\over{3}} \, \lambda^2 \,\, \equiv \,\, c_s \,  \lambda^2 \, . 
\monend  
%


\monitem 
At  second  order the equilibrium values have to be constrained: 
the ``heat flux'' at equilibrium is given by the relation 
\moneq   \label{d3q27-cns-1-ordre-2}   
q ^{\rm eq}  = 2 \, {{\sigma_5 - 4 \, \sigma_7 }\over{ \sigma_5 + 2 \  \sigma_7 }}   \, 
\lambda^2 \, J \,  
 \monend  
and a null value for the equilibrium   of the third order moments is imposed:
\moneq   \label{d3q27-cns-2-ordre-2}   
m_{23}^{\rm eq} \,= \, m_{24}^{\rm eq} \,= \, m_{25}^{\rm eq} \,=\, m_{26}^{\rm eq} 
\, = \, 0 \, . 
 \monend  
Moreover, the relations (\ref{d3q19-cns-2-ordre-2}) between   
H\'enon's parameters of second order moments   still have     
 to be imposed.
Then the second order equivalent equations are isotropic and the 
coefficients  of  (\ref{edp-ordre-2}) follow the non-traditional relations 
\moneq   \label{d3q27-coeffs-ordre-2}   
\mu \,=\, {{\sigma_5 \,\sigma_7}\over{\sigma_5 \,+ 2 \,\sigma_7}} \, \lambda \, \Delta x
\,, \quad \zeta_b \, \equiv \, 3 \, \zeta - \mu 
\,=\,  {{\lambda \, \sigma_4 \, \Delta x}\over{\sigma_5 \,+ 2 \,\sigma_7}} 
\, \big( \,  \sigma_5 -  2 \, \sigma_7  -  \alpha \, \sigma_5  - 2 \, \alpha \, \sigma_5   
\, \big) \, . 
 \monend  
%

\monitem 
At third order, we have   two options as for the D3Q19 scheme. 
If we suppose that   the heat flux at equilibrium  and 
only one time relaxation are  fixed, {\it id est}
\moneq   \label{d3q27-cn-1-ordre-3}  
\left\{\begin{array}{l}  \displaystyle 
q ^{\rm eq}  = -  2  \, \lambda^2 \, J  \,,  \quad 
 \varepsilon_2^{\rm eq}  \,=\,  - ( 2 + 3 \, \alpha )  \, \lambda^2 \, \rho  \, ,  
  \\  \vspace  {-.5 cm}   \\      \displaystyle 
 XX_e^{\rm eq} \, = \,  WW_e^{\rm eq}  \, = \, XY_e^{\rm eq} \, = \,  YZ_e^{\rm eq}  
\, = \,  ZX_e^{\rm eq}  \,=\, 0     \,, \quad 
\sigma_5 \,=\, \sigma_7 \,,   
   \end{array} \right. \monend  
the third order equivalent equations of the D3Q27 lattice Boltzmann scheme 
are   given by the expressions   (\ref{edp-ordre-3}). 
The coefficients in these equations are simple to evaluate 
with a software  of formal calculus:
\moneq   \label{d3q27-coefs-1-ordre-3}   
\left\{\begin{array}{l}  \displaystyle  
\mu     = {{ 1  }\over { 3 }}  \, \sigma_5 \,  \lambda \, \Delta  x  \, , \quad  
\zeta_b =  -{{1}\over{3}} \, \sigma_4 \, (1 + 3 \, \alpha)  \,  \lambda \, \Delta  x  \,, \quad 
\xi  \,=\, {{1}\over{36}} \, ( \alpha - 1 ) \, \Delta x^2 \,,    
  \\  \vspace  {-.5 cm}   \\      \displaystyle   \chi  = 
 {{1}\over{54}} \, ( \alpha + 2 ) \, \big( \,  1 +  \alpha +  6 \,\alpha \,\sigma_4^2   
  -2 \, \sigma_4^2 + 8 \,\sigma_5^2 \, \big)   \, \lambda^2 \,  \Delta x^2  \, . 
\end{array}\right. \monend   
%
 
\monitem 
The second solution for third order isotropy 
does not specify completely the ``square 
of the energy'' $ \, \varepsilon_2 \,$  
but fixes an important number of relaxation times: 
\moneq   \label{d3q27-cn-2-ordre-3} 
\left\{\begin{array}{l}  \displaystyle 
q ^{\rm eq}  =   -  2  \,  \lambda^2 \, J  \,,  \quad 
 \varepsilon_2^{\rm eq}  \,=\,  \beta   \, \lambda^2 \, \rho  \, ,   \quad 
\sigma_5 \,=\, \sigma_7 \,, 
  \\  \vspace  {-.5 cm}   \\      \displaystyle  
\sigma_{10} \,=\, \sigma_{11} \,=\, \sigma_{12} \,=\, {{1}\over{12 \, \sigma_5}} \,, \quad 
\sigma_{23} \,=\,  {{1}\over{12 \, \sigma_5}} \,, \quad 
\sigma_{24} \,=\,  {{1}\over{4 \, \sigma_5}} \,, \quad 
\sigma_{25} \,=\,  {{1}\over{8 \, \sigma_5}} \,, \quad 
  \\  \vspace  {-.5 cm}   \\      \displaystyle 
 XX_e^{\rm eq} \, = \,  WW_e^{\rm eq}  \, = \, XY_e^{\rm eq} \, = \,  YZ_e^{\rm eq}  
\, = \,  ZX_e^{\rm eq}  \,=\, 0   \,. 
   \end{array} \right. \monend   
The parameters $ \, \mu$,  $ \, \zeta$,  $ \, \zeta_b\,$ and   $ \, \xi \, $
are still given by the relations   (\ref{d3q27-coefs-1-ordre-3}). But the value 
of the parameter $ \, \chi \,$ is modified:  
\moneq   \label{d3q27-coefs-2-ordre-3}   
\left\{\begin{array}{l}  \displaystyle  
\chi  = 
 {{1}\over{162 \, \sigma_5 }} \, \Big( \,  4 \, \sigma_5  + 2 \, \sigma_4  + 3 \, \alpha \, \sigma_4
 - 6 \, \alpha \, \sigma_5   - 3 \, \alpha^2 \, \sigma_5
+ 18  \, \alpha^2 \, \sigma_4^2 \, \sigma_5   + 42 \, \alpha\, \sigma_4^2   \, \sigma_5 
  \\  \vspace  {-.5 cm}   \\      \displaystyle \qquad  \qquad  \qquad 
           + 12 \, \sigma_4^2 \, \sigma_5      
   - 24 \, \alpha \, \sigma_5^3  - 48 \, \sigma_5^3    
  + \beta\, \sigma_4 - \beta \, \sigma_5\, \Big)   \, \lambda^2 \,  \Delta x^2  \, . 
\end{array}\right. \monend

\monitem 
The search of an isotropic form of the fourth order equivalent      
partial differential equations like   (\ref{edp-ordre-4})  
leads to a nonlinear system of 33 equations.        
We have obtained a first solution 
 with the following  particular parameters: 
\moneq   \label{d3q27-cn-1-ordre-4}  
  r ^{\rm eq}  = 2 \, \lambda^4 \, J \,, \quad 
 \sigma_{10}  =   \sigma_{11}  =   \sigma_{12}  \,, \quad 
  \sigma_{18}  =   \sigma_{19}  \,, \quad 
 \sigma_{23}  =   \sigma_{24}  =   \sigma_{25}  \, . 
 \monend   
It is possible   to fix the other parameters of the scheme
with the ratio $ \, \psi \,$ of   H\'enon's parameters 
associated with the moments $ \, XX_e \, $ and $\, XX$. 
With the notations proposed in  (\ref{d3q27-moments}), 
we set  
\moneq   \label{d3q27-psi-ordre-4}  
\psi \, \equiv \, {{\sigma_{18}}\over{\sigma_5}} \, . 
 \monend   
When we impose  the following relations between the coefficients of relaxations
(all defined through their associated H\'enon's parameter
introduced at the relation (\ref{d2q9-cns-ordre-2})),  
\moneq   \label{d3q27-cn-2-ordre-4} 
\left\{\begin{array}{l}  \displaystyle 
\sigma_4 \,= \, \sigma_5 \, \, 
{{ 3 \, \psi^3 - 4 \, \psi^2   - 13 \, \psi + 32 }\over{ 3 \, \psi^3 - 22 \, \psi^2   + 23 \, \psi + 14   }} 
\, ,   \quad  
  \\  \vspace  {-.5 cm}   \\      \displaystyle 
\sigma_{10} \,= \, \sigma_{11} \,=  \, \sigma_{12} \,=  \, {{1}\over{12 \, \sigma_5 }}
\, \,  {{ (3 \, \psi - 7) \, ( \psi - 4) }\over{ 3 \, \psi^2 - 11 \, \psi + 14   }} \,, 
  \\  \vspace  {-.5 cm}   \\      \displaystyle
 \sigma_{16} \,= \,  \sigma_5 \, \, 
{{ 6 \, \psi^5 - 24 \, \psi^4 + 100 \, \psi^3 - 267  \, \psi^2 + 506 \, \psi - 364  }\over{ 
(4 \, \psi - 7) \, ( 3 \, \psi^3 - 22 \, \psi^2   + 23 \, \psi + 14 )  }} \,, 
  \\  \vspace  {-.5 cm}   \\      \displaystyle
 \sigma_{18} \,= \,  \sigma_{19} \, \equiv  \, \psi \,  \sigma_5  \,, 
  \\  \vspace  {-.5 cm}   \\      \displaystyle
 \sigma_{23} \,= \,  \sigma_{24}   \,= \,  \sigma_{25} \,=\,  {{1}\over{12 \, \sigma_5  }}
\,\, {{ 3 \, \psi^2 - 7 \, \psi + 16  }\over{ 3 \, \psi^2 - 11 \, \psi + 14   }} \,, 
  \\  \vspace  {-.5 cm}   \\      \displaystyle
 \sigma_{26} \,= \,  {{1}\over{18 \, \sigma_5  }}
\,\,  {{ 21 \, \psi^2 - 73 \, \psi + 100  }\over{ 3 \, \psi^2 - 11 \, \psi + 14   }}   \,,  
   \end{array} \right. \monend   
the equivalent partial differential equations  of  the D3Q27 lattice Boltzmann scheme 
are isotropic at fourth order of accuracy. We can make more explicit   
graphically the previous     
result.
We observe in Figure~1 that the fundamental stability property $ \, \sigma_{16} > 0 \,$ can   be     
maintained only if  $ \,  0  <  \psi  <  1.5  $. 
With this restriction, we see in Figure~1     
again that   H\'enon's parameters 
 $ \, \sigma_{10}, $  $ \, \sigma_{16}  \,$ and  $ \, \sigma_{18}  \,$ 
remain positive    only if 
\moneq   \label{d3q27-cn-3-ordre-4}   
 1 \,< \,   {{ \sigma_{4} }\over{ \sigma_{5} }} \, < \, 2.25  \, . 
\monend    
%
\smallskip   \smallskip                 
\begin {center}
\begin{picture}(0,0)%
\includegraphics{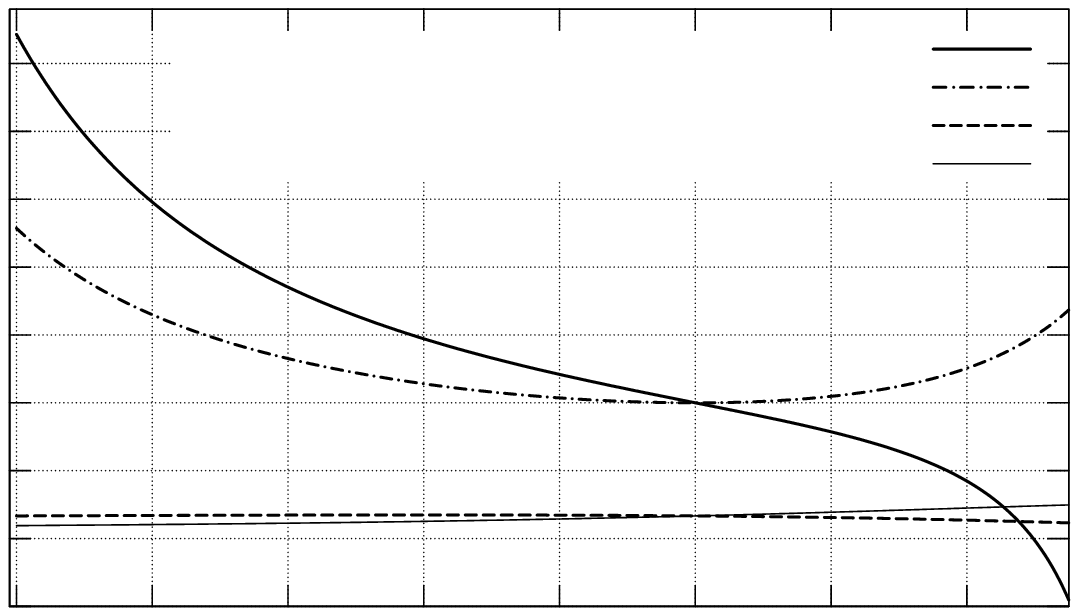}%
\end{picture}%
\begingroup
\setlength{\unitlength}{0.0200bp}%
\begin{picture}(18000,10800)(0,0)%
\put(1650,1650){\makebox(0,0)[r]{\strut{}-0.5}}%
\put(1650,2627){\makebox(0,0)[r]{\strut{} 0}}%
\put(1650,3605){\makebox(0,0)[r]{\strut{} 0.5}}%
\put(1650,4582){\makebox(0,0)[r]{\strut{} 1}}%
\put(1650,5559){\makebox(0,0)[r]{\strut{} 1.5}}%
\put(1650,6536){\makebox(0,0)[r]{\strut{} 2}}%
\put(1650,7514){\makebox(0,0)[r]{\strut{} 2.5}}%
\put(1650,8491){\makebox(0,0)[r]{\strut{} 3}}%
\put(1650,9468){\makebox(0,0)[r]{\strut{} 3.5}}%
\put(2023,1100){\makebox(0,0){\strut{} 0}}%
\put(3978,1100){\makebox(0,0){\strut{} 0.2}}%
\put(5933,1100){\makebox(0,0){\strut{} 0.4}}%
\put(7888,1100){\makebox(0,0){\strut{} 0.6}}%
\put(9843,1100){\makebox(0,0){\strut{} 0.8}}%
\put(11798,1100){\makebox(0,0){\strut{} 1}}%
\put(13754,1100){\makebox(0,0){\strut{} 1.2}}%
\put(15709,1100){\makebox(0,0){\strut{} 1.4}}%
\put(9550,275){\makebox(0,0){\strut{} $\psi \equiv \sigma_{18} \, / \, \sigma_5$  }}%
\put(14950,9675){\makebox(0,0)[r]{\strut{}$\sigma_{16}  \, / \, \sigma_5   \,\,$ }}%
\put(14950,9125){\makebox(0,0)[r]{\strut{}$\sigma_4 \, / \, \sigma_5 \,\,$ }}%
\put(14950,8575){\makebox(0,0)[r]{\strut{}$\sigma_{10}  \, / \, \sigma_5   \,\,$ }}%
\put(14950,8025){\makebox(0,0)[r]{\strut{}$\sigma_{26} \,  / \, \sigma_5   \,\,$ }}%
\end{picture}%
\endgroup

\end {center}
\smallskip \noindent  {\bf Figure 1}. \quad  Fourth order isotropy parameters for the 
D3Q27 lattice Boltzmann scheme. 
\smallskip \smallskip 

%
\noindent
Due to the expressions (\ref{d3q27-coefs-1-ordre-3}) of the shear and the bulk
viscosities, the inequality (\ref{d3q27-cn-3-ordre-4}) imposes signifiant   
 restrictions for the physical parameters $ \, \mu \,$ and $ \, \zeta_b .$ 
The coefficients $ \, \eta \,$ and $ \, \mu_4 \,$ associated 
with the fourth order equation    
  (\ref{edp-ordre-4})  can be evaluated easily: 
\moneq   \label{d3q27-coeff-1-ordre-4}  
\left\{    \begin{array}{l}  \displaystyle  
\eta  \,\,\,\, = \,   {{\sigma_5\, \lambda  \, \Delta x^3}\over{108}}     {{  N_\eta  }\over
{ 14 + 23\, \psi - 22 \, \psi^2 + 3 \, \psi^3 }}  , 
  \\  \vspace  {-.5 cm}   \\      \displaystyle
N_\eta \, = \,  (\alpha+2) \, (
32 \, \alpha - 8 - 13 \, \alpha \, \psi   - 35 \, \psi 
  - 4 \, \alpha \, \psi^2   + 28 \, \psi^2 + 3 \, \alpha \, \psi^3 - 3 \, \psi^3 )
  \\  \vspace  {-.5 cm}   \\      \displaystyle
 \mu_4  \,\, = \,  {{\sigma_5 \, \lambda \, \Delta x^3}\over{108}} \, 
{{ 132 \, \psi \, \sigma_5^2 - \psi  + 36 \,  \psi^2 \, \sigma_5^2 + 168  \, \sigma_6^2 
+ 3 \, \psi^2 - 8   }\over{ 3 \, \psi^2 - 11 \, \psi + 14   }} \, .
   \end{array} \right. \monend   
The expression of $\, \zeta_4 \,$ is quite long and is       
reported in the relation (\ref{d2q27-zeta4-ordre-4}) of            
 Appendix~2. 

  \bigskip \bigskip   \noindent {\bf \large 8) \quad  Conclusion   } 

\monitem  
In this contribution, we have presented the                                  
``Berliner version'' of the   Taylor expansion method in the linear case. 
This is done with explicit algebra and allows a 
huge reduction of  computer time for formal analysis.      
We have also considered in all generality 
 acoustic type partial differential equations that are 
rotationally invariant at an arbitrary order.             

\monitem  
The generalization of a methodology of group theory 
 for discrete invariance groups
 of a lattice Boltzmann scheme
remains still under question, in the spirit of the previous
study of  Rubinstein and  Luo \cite{RL08}.


\monitem  
Concerning the fundamental examples considered in this contribution,      
the   D2Q9  scheme   can be   invariant by rotation at   third order.
At fourth  order, physical parameters have to be strongly correlated.  
The    D2Q13   scheme is invariant by rotation  at    fourth order 
for an {\it ad hoc} fitting of the parameters.  
We have not explored all the possible solutions of the strongly nonlinear set of equations
that is necessary to solve in order to fit the fourth order isotropy.      
Numerical experiments have to confirm our theoretical considerations. 
The D3Q19  lattice Boltzmann scheme admits two sets  of 
coefficients in order to impose rotational  invariance  
at third order. Particular physics has
to be imposed to satisfy  fourth order isotropy.  
The D3Q27 scheme  is  rotationally invariant at fourth order        
for a parameterized set of parameters. Our analysis imposes 
restrictions for the physical parameters to   guarantee the stability. 
A complementary numerical experiment will be welcome !        

\bigskip  \bigskip   \noindent {\bf \large Acknowledgments}   

 \noindent     
This work has been financially supported by the French Ministry of Industry (DGCIS) 
and the ``R\'egion Ile-de-France''   in the framework of the LaBS Project. 
The authors thank  Yves Benoist (Centre National de la Recherche Scientifique
and department of Mathematics in Orsay) for an enlighting  discussion
about the representation of groups.  
The authors thank  also the referees for helpful comments and suggestions.  

\bigskip \bigskip     \newpage 
 \noindent {\bf \large  Appendix 1. Formal expansion in the linear case } 

\monitem 
We present in this Appendix the  ``Berliner version'' \cite{Du11}   of the algorithm 
proposed in all generality in our contribution \cite{DL09}. 
We suppose having defined a 
lattice Boltzmann scheme  ``DdQq'' with  $d$ space dimensions 
and  $q$ discrete velocities at each vertex. 
The invertible matrix  $M$ between the particules and the moments is given:
\moneq   \label{ap1-matrice-M}  
m_k \,=\,  \sum_{j=0}^{q-1}   M _{kj} \, f_j  \,\, \equiv \,\, \big( M \smb f \big)_k \,,
\quad 0 \leq k \leq q - 1 \, . 
\monend 
The lattice Boltzmann scheme generates 
$N$  conservation laws: the  first moments 
\moneqstar
 m_k  \, \equiv \, W_k \,, \qquad  0 \leq k  \leq N-1 
\monendstar
are conserved  during the collision step~:  
\moneq   \label{ap1-moments-conserves}  
m^*  \, = \, m_k \, = \,  W_k \,  . 
\monend 
The  $q-N$ ``slave'' moments  $\, Y \, $ with
\moneq   \label{ap1-moments-esclaves}  
 Y_\ell \, \equiv \, m_{N+\ell} \,, \qquad  0 \leq \ell \leq q-N-1
\monend 
relax   towards an equilibrium value   $\, Y_\ell^{\rm eq}  . \,$ 
This   equilibrium value  
is supposed to be a {\bf  linear }  function of the state $W$.
We introduce a constant rectangular matrix $\, E \, $ with $ \, N-q\,$ 
lines and $\, N \,$ columns to represent this linear function:   
\moneq   \label{ap1-linearite}  
 Y_\ell^{\rm eq} \,=\,  \sum_{k=0}^{N-1} 
 E_{\ell k} \,\, W_k \,, \qquad  0 \leq \ell \leq q-N-1 \, .
\monend 
The relaxation step is obtained through the usual algorithm \cite{DdH92}
that decouples the  moments:   
\moneq   \label{ap1-relaxation-Y}  
  Y_\ell^*  \,=\,   Y_\ell + s_\ell \, ( Y_\ell^{\rm eq} -   Y_\ell ) \,, 
\,\, s_\ell > 0 \,, \qquad  0 \leq \ell \leq q-N-1 \, .
\monend 
Observe that the numbering of the ``$s$'' coefficients used in (\ref{ap1-relaxation-Y})
differ  just a little from the one used for the equation (\ref{relaxation}) and 
the four examples considered previously. 
With a matricial notation, the relaxation can be written as:  
\moneq   \label{ap1-relaxation-m}  
   m^* \,=\,  J_0   \, \smb\, m 
\monend 
with a matrix  $J_0$ of order  $q$ decomposed by blocks according to 
\moneq   \label{ap1-matrice-J0}  
J_0    \,=\,  \begin {pmatrix} \displaystyle  {\rm I}_N  &  0   
 \\[3mm]   \displaystyle  S \, \smb\, E   &  \,\,  {\rm I}_{q-N} - S  
\end  {pmatrix}  \,  
\monend 
and    a diagonal matrix $ \, S \,$ of order $\, q-N \,$  defined by 
$\,  S \,\equiv\, {\rm diag}  \big(s_0 \,,\, s_1 \,,\, \dots \,, s_{q-N-1} \big ) .\, $ 
The discrete advection step  follows the method of characteristics: 
\moneq   \label{ap1-evolution-temps} 
f_j(x, t+ \Delta t) \,=\,   f_j^*(x- v_j \Delta t, t) \, , \qquad 0 \leq j \leq q-1 \, .  
\monend 

\monitem
With the d'Humi\`eres's lattice Boltzmann scheme \cite{DdH92} 
previously defined, we can 
proceed to a formal Taylor expansion: 
 
 \smallskip    \noindent 
$ \displaystyle m_k(t+ \Delta t)  \,=\,   
\sum_{j} M_{k j} \,  f_j^*(x- v_j \Delta t)  \,=\, 
\sum_{j \ell} M_{k j} \,  M^{-1}_{j \ell} \,  m_{\ell}^*(x- v_j \Delta t)    $ 

$ \displaystyle \qquad \qquad   \,=\, \sum_{j \ell} M_{k j} \,  M^{-1}_{j \ell} \, 
 \sum_{n=0}^{\infty}
 {{\Delta t^n}\over{n \, !}} \, \bigg( - \sum_{\alpha=1}^d    v_j^{\alpha} \, \partial_{\alpha} 
\bigg)^n  \, m_{\ell}^*   $ 

$ \displaystyle \qquad \qquad   \,=\,   \sum_{n=0}^{\infty}  {{\Delta t^n}\over{n \, !}} \,
 \sum_{j \ell p } M_{k j} \,  M^{-1}_{j \ell  } \, 
  \bigg( - \sum_{\alpha=1}^d    v_j^{\alpha} \, \partial_{\alpha} 
\bigg)^n  \,( J_0 )_{\ell p} \, m_{p}   \, .  $ 

 \smallskip    \noindent 
We introduce a derivation matrix of order  $n \geq 0$, defined by blocks 
of space differential operators of order  $n$: 
\moneq   \label{ap1-matrice-derivation} 
 \begin {pmatrix} \displaystyle  A_n  &  B_n    
 \\[3mm]   \displaystyle  C_n  &  D_n   
\end  {pmatrix}_{k \, p}   \, \equiv  \,  {{1}\over{n \, !}}  \,  
\sum_{  j  \, \ell} M_{k \, j} \,  \big( M^{-1} \big)_{j \, \ell } \,\, 
  \bigg( - \sum_{\alpha=1}^d    v_j^{\alpha} \, \partial_{\alpha} 
\bigg)^n  \,( J_0 )_{\ell p} \,, \quad n \geq 0 \, .  
\monend
We observe that in the relation (\ref{ap1-matrice-derivation}), the blocks $\, A_n \,$ 
and $\, D_n \,$ are square matrices of order $N$ and $ \, q-N \,$ respectively.
The matrices  $\, B_n \,$ and $\, C_n \,$ are rectangular of order  
$ \,  N \times ( q-N )  \,$  and $ \,  ( q-N ) \times N \,$  
respectively. We remark also that at order zero, the matrices       
 $\, A_0 , \,$  $\, B_0 , \,$  $\, C_0  \,$ and   $\, D_0  \,$ 
are known: 
\moneq   \label{ap1-bloc-zero} 
 \begin {pmatrix} \displaystyle  A_0  &  B_0    
 \\[3mm]   \displaystyle  C_0  &  D_0   \end  {pmatrix}  \, = \,
  J_0   \,=\,  \begin {pmatrix} \displaystyle  {\rm I}_N  &  0   
 \\[3mm]   \displaystyle  S \, \smb\, E   &  \,\,  {\rm I}_{q-N} - S  
\end  {pmatrix}  \,  .
\monend
The previous Taylor expansion can now be written under a matricial form:
\moneq   \label{ap1-Taylor-matrice} 
  \begin {pmatrix} \displaystyle W 
 \\[3mm]   \displaystyle  Y \end  {pmatrix} (x, \, t+\Delta t)  = 
 \sum_{n=0}^{\infty}  \, \Delta t ^n \,  
 \begin {pmatrix} \displaystyle  A_n  &  B_n    
 \\[3mm]   \displaystyle  C_n  &  D_n   \end  {pmatrix}  
 \, \smb \,  \begin {pmatrix} \displaystyle W 
 \\[3mm]   \displaystyle  Y \end  {pmatrix} (x, \, t )  \, . 
\monend

\monitem 
At  order zero  relative to  $ \, \Delta t$  we have:     
\moneqstar   \begin {pmatrix} \displaystyle W  \\
 \displaystyle  Y \end  {pmatrix} (x, \, t) \,+\, {\rm O}(\Delta t) \,=\, 
J_0  \, \smb \,  \begin {pmatrix} \displaystyle W  \\
  \displaystyle  Y \end  {pmatrix}  \,+\, {\rm O}(\Delta t)   \,=\, 
 \begin {pmatrix} \displaystyle W 
 \\    \displaystyle S  \, \smb \, E \, \smb \,   W 
+ ( {\rm I} -S ) \,\smb \, Y  \end {pmatrix}     
   \,+\, {\rm O}(\Delta t) 
\monendstar  
  and   the non-conserved moments are close to the equilibrium:
\moneq   \label{ap1-proche-equilibre} 
 Y(x,\, t) \,= \,  E \, \smb \,   W (x,\, t) \,+\,  {\rm O}(\Delta t) \, . 
\monend

\monitem We make now the hypothesis of a 
general form for the expansion of the nonconserved moments:   
\moneq   \label{ap1-developpement-Y} 
Y(x,\, t) \,= \,  \Big( E + \sum_{n \geq 1} \Delta t^n \, \beta_n \Big)  
 \, \smb \,   W (x,\, t) 
\monend 
and the hypothesis of a formal linear partial differential system 
of arbitrary order for the conserved variables $W$:   
\moneq   \label{ap1-edp-W} 
  {{\partial W}\over{\partial t}} \,=\,    \Big( \sum_{\ell  \, \geq \, 0} \Delta t^\ell \, 
 \alpha_{\ell+1}   \Big)   \, \smb \,   W (x,\, t)  \,,
\monend 
where   $\alpha_{\ell}$ and  $\beta_n$ are space differential operators 
of order $\ell$ and $n$ respectively. 
We develop the first equation of  (\ref{ap1-Taylor-matrice}) up to first  order:  

\smallskip   \noindent \qquad  \qquad   $ \displaystyle
W \,+\,  \Delta t \, {{\partial W}\over{\partial t}}   \,+\,   {\rm O}(\Delta t^2)  \,=\,
W \,+\,  \Delta t \, \big( A_1 \, W + B_1 \, Y \big)   \,+\,   {\rm O}(\Delta t^2) $ 

\smallskip   \noindent \qquad  \qquad  
$ \displaystyle  \qquad  \qquad  \qquad  \qquad  \qquad  \quad   \,=\,
W \,+\,  \Delta t \, \big( A_1 \, W + B_1 \, E \, W  \big)   \,+\,   {\rm O}(\Delta t^2) $ 

\smallskip   \noindent
 due to (\ref{ap1-proche-equilibre}). 
Then    
\moneq   \label{ap1-edp-ordre-1} 
 {{\partial W}\over{\partial t}}  \,=\, 
\big( A_1   + B_1 \, E    \big)  \,\smb\,  W \,+\,   {\rm O}(\Delta t) 
\monend 
and the relation (\ref{ap1-edp-W}) 
is satisfied at order  one, with
\moneq   \label{ap1-alpha-1} 
 \alpha_1 \,=\, A_1 \,+\, B_1 \, E \, . 
\monend 
The  ``Euler equations''   are emerging !  
We have an analogous calculus 
for the second equation of  (\ref{ap1-Taylor-matrice}) : 
\moneqstar 
Y    \,+\,  \Delta t \, {{\partial Y}\over{\partial t}}    \,+\,   {\rm O}(\Delta t^2)
\,=\, 
S \, E \, W   \,+\,    ({\rm I} - S) \, Y \,+\, \Delta t \,
  \big( C_1 \, W  + D_1 \, E \, W  \big)     \,+\,   {\rm O}(\Delta t^2) \, . 
\monendstar 
  We clarify    
  the time derivative $ \,  \partial_t  Y   \,$ 
at order zero   by differentiating (formally !)  the relation  (\ref{ap1-proche-equilibre}) 
relative to time:                  
\moneqstar 
 {{\partial Y}\over{\partial t}} \,=\,  E \,   {{\partial W}\over{\partial t}}  
\,+\,   {\rm O}(\Delta t) \,= \, E \, \alpha_1 \, W \,+\,   {\rm O}(\Delta t)   \, . 
\monendstar 
We introduce this expression inside the previous calculus. Then:  
\moneqstar 
S \, Y  \,+\, \Delta t \,  E \, \alpha_1 \, W \,+\,   {\rm O}(\Delta t^2) \,=\, 
S \, E \, W  \,+\, \Delta t \,  \big( C_1 \, W  + D_1 \, E \, W  \big) 
    \,+\,   {\rm O}(\Delta t^2) \, . 
\monendstar 
Consequently we have established the expansion of the nonconserved moments at order one: 
\moneq   \label{ap1-Y-ordre-1} 
 Y  \,=\, E \, W \,+\, \Delta t \, S^{-1} \, 
 \big( C_1   + D_1 \, E - E \, \alpha_1 \big) \, W  \,+\,   {\rm O}(\Delta t^2) 
\monend
with 
\moneq   \label{ap1-beta-1} 
 \beta_1 \,=\, S^{-1} \,  \big( C_1   + D_1 \, E - E \, \alpha_1 \big) \,.
\monend
Now, we have formally  

\smallskip \noindent     $ \displaystyle  
   {{\partial^2 W}\over{\partial t^2 }}  \,=\, 
 {{\partial}\over{\partial t}} \big( \alpha_1 \, W  \,+\,   {\rm O}(\Delta t) \big)  \,=\,
\alpha_1 \,  {{\partial W}\over{\partial t}} \,+\,   {\rm O}(\Delta t) 
 \,=\, \alpha_1 \,  \big( \alpha_1 \, W \big) \,+\,   {\rm O}(\Delta t) \,=\,
     \alpha_1^2 \, W  \,+\,   {\rm O}(\Delta t) $ 

\smallskip \noindent  
and we recognize the ``wave equation''   
\moneq   \label{ap1-wave} 
   {{\partial^2 W}\over{\partial t^2 }}  \,-\,   \alpha_1^2 \, W  \,=\,   {\rm O}(\Delta t) \, . 
\monend
%

\monitem  We can derive a formal expansion at order two.
We go one step further in the Taylor expansion of equation   (\ref{ap1-Taylor-matrice})   : 

\smallskip \noindent    $ \displaystyle
W \,+\,  \Delta t \, {{\partial W}\over{\partial t}}   \,+\, 
{1\over2} \,   \Delta t^2  \, \alpha_1^2 \, W \,+\,   {\rm O}(\Delta t^3)  \,=\, $ 

\smallskip   $ \displaystyle  \qquad \qquad \,=\,
W \,+\,  \Delta t \, \big( A_1 \, W + B_1 \, Y \big)   \,+\,  
 \Delta t^2  \, \big( A_2 \, W + B_2 \, Y \big)  \,+\,    {\rm O}(\Delta t^3) $  

\smallskip   $ \displaystyle  \qquad \qquad \,=\,
W \,+\,  \Delta t \, \big( A_1 \, W + B_1 \, (E \, W +  \Delta t \, \beta_1 \, W)  \big)   \,+\,  
 \Delta t^2  \, \big( A_2 \, W + B_2 \, E \, W \big)  \,+\,    {\rm O}(\Delta t^3) $ 

\smallskip \noindent 
and dividing by  $ \, \Delta t $, we obtain a 
 ``Navier-Stokes  type'' second order equivalent equation: 
\moneqstar  
 {{\partial W}\over{\partial t}}  \,=\, 
\alpha_1 \, W \, + \,   \Delta t \, \bigg(  B_1 \,  \beta_1 + A_2 +  B_2 \, E
- {1\over2} \,  \alpha_1^2 \bigg) \, W  \,+\,    {\rm O}(\Delta t^2) \, .
\monendstar
With the notations introduced in  (\ref{ap1-edp-W}),
 we have made explicit         
the partial differential equations for the conserved variables at the order two: 
\moneqstar  
 {{\partial W}\over{\partial t}}  \,=\, 
\alpha_1 \, W \, + \,   \Delta t \, \alpha_2 \,  W  \,+\,    {\rm O}(\Delta t^2) 
\monendstar
with 
\moneq   \label{ap1-alpha-2} 
 \alpha_2  \,=\,    A_2 +  B_2 \, E \,+\,  B_1 \,  \beta_1 \,-\, {1\over2} \,  \alpha_1^2  \, . 
\monend 
We remark that this Taylor expansion method can be viewed as 
a  ``numerical  Chapman Enskog expansion'' relative to a specific 
numerical parameter $ \, \Delta t \, $ instead of a small physical relaxation time step. 
For the moments $ \, Y \,$ out of equilibrium,   
we expand the first order derivative of  $Y$ relative to time 
with a formal derivation of the relation (\ref{ap1-Y-ordre-1}):

\smallskip  \noindent $ \displaystyle \qquad \qquad 
 {{\partial Y}\over{\partial t}}  \,=\, 
 {{\partial}\over{\partial t}} \Big( E \, W \,+\, 
\Delta t  \, \beta_1 \, W \Big) \,+  {\rm O}(\Delta t^2) $  

\smallskip  \noindent $ \displaystyle   \qquad \qquad  \qquad  \,=\, 
 E \, \big( \alpha_1 \, W \, +\, \Delta t  \,  \alpha_2  \, W \big) 
\, +\, \Delta t  \, \beta_1 \,  \alpha_1 \, W \,+\,   {\rm O}(\Delta t^2) $ 

\smallskip  \noindent $ \displaystyle   \qquad \qquad  \qquad  \,=\, 
 \Big(  E \,  \alpha_1 \,+\, \Delta t \, 
\big(  E \,  \alpha_2   \,+\, \beta_1 \, \alpha_1 \big) \, \Big) \, W 
 \,+\,   {\rm O}(\Delta t^2) \, . $ 

 \noindent Then 
\moneq   \label{ap1-dtY-ordre-1} 
 {{\partial Y}\over{\partial t}}  \,=\, 
 \Big(  E \,  \alpha_1 \,+\, \Delta t \, 
\big(  E \,  \alpha_2   \,+\, \beta_1 \, \alpha_1 \big) \, \Big) \, W 
 \,+\,   {\rm O}(\Delta t^2) \, . 
\monend 

\smallskip  \noindent
Analogously for the second order time derivative: 
%
\moneq   \label{ap1-dt2Y-ordre-0}  
{{\partial^2 Y}\over{\partial t^2}}  \,=\,
 E \, \alpha_1^2 \, W   \,+\,   {\rm O}(\Delta t)   \, .   
\monend 
We re-write the second line  
 of the expansion of the equation   (\ref{ap1-Taylor-matrice}) at second order accuracy:

\smallskip  \noindent  $ \displaystyle 
Y  \,+\, \Delta t \,  {{\partial Y}\over{\partial t}} \,+\, {{ \Delta t^2}\over{2}} 
 \, {{\partial^2 Y}\over{\partial t^2}}   \,+\,   {\rm O}(\Delta t^3) \,=\,  $ 

\smallskip  \hfill   $ \displaystyle  = \, 
S \, E \, W \,   \,+\, ({\rm I} - S) \, Y 
\,+\,  \Delta t \, \big( C_1 \, W \,+\, D_1 \, Y) 
\,+\,  \Delta t^2 \, \big( C_2 \, W \,+\, D_2 \, Y)   \,+\,   {\rm O}(\Delta t^3)   $  

\noindent
and we get 

\smallskip   \noindent $ \displaystyle  S  Y \,=\, 
S \, E \, W \,-\,\Delta t \,   \big(  E \,  \alpha_1 \,+\, \Delta t \, 
(  E \,  \alpha_2   \,+\, \beta_1 \, \alpha_1  ) \, \big) \, W 
\,-\, {{ \Delta t^2}\over{2}} \, E \, \alpha_1^2 \, W     $

\smallskip   \noindent \qquad  \quad   $ \displaystyle 
\,+\,  \Delta t \, \Big( C_1 \, W \,+\, D_1 \, \big( E + \Delta t \, \beta_1 \big) \, W \Big) 
\,+\,  \Delta t^2 \, \big( C_2 \, W \,+\, D_2 \, E \, W)   \,+\,   {\rm O}(\Delta t^3)   $ 
 
\smallskip   \noindent  $ \displaystyle   Y \,\,\,\, \,=\,  E \, W \,+\, \Delta t  \,\, S^{-1} \, 
\big( C_1 \,+\, D_1 \, E \,-\,  E \, \alpha_1 \big) \, W  $ 
 
\smallskip   \noindent \qquad  \quad   $ \displaystyle 
\,+\,   \Delta t^2  \,\, S^{-1} \, 
\Big( C_2 \,+\, D_2 \, E \,  +\, D_1 \, \beta_1 \,-\, E \, \alpha_2 \,-\, \ \beta_1
\,\alpha_1 \,-\, {{1}\over{2}} \, E \, \alpha_1^2 \Big) \, W    \,+\,   {\rm O}(\Delta
t^3) \, .  $ 

\smallskip   \noindent  
It is exactly the expansion   (\ref{ap1-edp-W})  at second order :  
\moneqstar  
Y \,=\ E \, W \,+\,  \Delta t \, \beta_1  \,  W  \,+\,  \Delta t^2 \, \beta_2 \,  W 
 \,+\,    {\rm O}(\Delta t^2) 
\monendstar
with  
\moneq   \label{ap1-beta-2} 
 \beta_2  \,=\, S^{-1} \, \Big[ \, 
  C_2 \,+\, D_2 \, E \,  +\, D_1 \, \beta_1 \,-\, E \, \alpha_2 \,-\, \ \beta_1
\,\alpha_1 \,-\, {{1}\over{2}} \, E \, \alpha_1^2  \, \Big] 
\monend

\monitem  
For the  general case, we proceed by induction.
We suppose that the developments   (\ref{ap1-developpement-Y})  and  (\ref{ap1-edp-W})  
are correct up to the order  $k$, that is: 
\moneq   \label{ap1-hyp-recur}  
\left\{ \begin{array} [c]{rcl}   \displaystyle 
{{\partial W}\over{\partial t}}  &=&  \Big( \alpha_1 
\,+\,   \Delta t \, \alpha_2 \,+ \dots \,   \Delta t^{k-1} \, \alpha_k \Big) \, W 
 \,+\,    {\rm O}(\Delta t^k )    
 \\ \displaystyle   \vspace{-.5cm}  ~  \\ \displaystyle
Y  &=&   \Big( E \, +\,   \Delta t \, \beta_1 
\,+\,   \Delta t^2  \, \beta_2 \,+ \dots \,   \Delta t^{k} \, \beta_k \Big) \, W 
 \,+\,    {\rm O}(\Delta t^{k+1} ) \, .  
\end{array} \right. \monend  
We expand the relation  (\ref{ap1-Taylor-matrice}) at order $k+2$, 
we eliminate the zeroth order term and divide by  $ \, \Delta t$. We obtain  
\moneq   \label{ap1-calcul-recur}     
 {{\partial W}\over{\partial t}}  \,+\, \sum_{j=2}^{k+1} {{ \Delta t^{j-1}}\over{j !}}
\, \big( \partial_t^j W \big)  \,+\,    {\rm O}(\Delta t^{k+1} ) \,=\, 
 \sum_{j=1}^{k+1} \Delta t^{j-1}  \, \big( A_j \, W \,+\, B_j \, Y \big) 
 \,+\,    {\rm O}(\Delta t^{k+1}  ) \, . 
\monend  
%
%
The term  $ \,  \partial_t^j W \, = \, 
\big(  \sum_{\ell=1}^{\infty}   \Delta t^{\ell-1} \, \alpha_\ell \big) ^j \,\,$ 
on the left hand side of  (\ref{ap1-calcul-recur})  
 can be evaluated by taking the  formal power of 
the equation   (\ref{ap1-edp-W})  at the order  $j$. We define the coefficients 
$ \, \Gamma_m^j \,$ according to: 
\moneq   \label{ap1-serie-Gamma}   
\Big(  \sum_{\ell=1}^{\infty}   \Delta t^{\ell-1} \, \alpha_\ell \Big) ^j \,\equiv \, 
 \sum_{\ell=0}^{\infty}   \Delta t^{\ell} \,\, \Gamma^j_{j+\ell} 
\, , \qquad j \geq 0 \,  .  
\monend  
They can be evaluated without difficulty from 
the coefficients $ \, \alpha_\ell$, taking care of the non-commutativity of 
the product of two matrices. 
We report the corresponding terms and we  identify  the coefficients
in factor of   $ \, \Delta t^{k} \,$  between the two sides
of the equation  (\ref{ap1-calcul-recur}),
with the help of the induction hypothesis
 (\ref{ap1-hyp-recur}). We deduce:   
\moneq   \label{ap1-recurrence-alpha}    
\alpha_{k+1} \,=\, A_{k+1}  \,+\,  \sum_{j=1}^{k+1} B_j \,\,  \beta_{k+1-j} 
 \,-\,   \sum_{j=2}^{k+1}  {{1}\over{j !}} \,  \Gamma^j_{k+1} \,. 
\monend  
We do the same operation   with the second relation of  (\ref{ap1-Taylor-matrice}) : 
\moneq   \label{ap1-dvpmt-Y}   
\left\{ \begin{array} [c]{c}   \displaystyle 
Y \,+\,   \sum_{j=1}^{k+1} {{ \Delta t^{j}}\over{j !}}
\, \big( \partial_t^j Y \big)  \,+\,    {\rm O}(\Delta t^{k+2} ) \,=\, \hfill 
 \\ \displaystyle   \vspace{-.5cm}  ~  \\ \displaystyle  
\qquad \qquad    \,=\,
 S\, E \, W \,+\, ({\rm I}-S) \, Y   \,+\, 
 \sum_{j=1}^{k+1} \Delta t^{j}  \, \big( C_j \, W \,+\, D_j \, Y \big) 
 \,+\,    {\rm O}(\Delta t^{k+2}  )  \, . 
\end{array} \right. \monend  
As in the previous case, we suppose that we have evaluated       
   formally the temporal derivative  

\smallskip  \noindent      $ \displaystyle \partial_t^j Y \,=\,  
\partial_t^j  \, \Big[ \big(  E \, +\,   \Delta t \, \beta_1 
\,+\,   \Delta t^2  \, \beta_2 \,+ \dots \,   + \Delta t^{k} \, \beta_k + \dots \big) \, W \Big]  $ 

\smallskip  \noindent      $ \displaystyle \qquad 
 \,=\, \big(  E \, +\,   \Delta t \, \beta_1 
\,+\,   \Delta t^2  \, \beta_2 \,+ \dots \,   + \Delta t^{k} \, \beta_k + \dots \big) \, 
\big(   \partial_t^j  W \big) \,$ 

\smallskip  \noindent      $ \displaystyle \qquad 
 \,=\, \big(  E \, +\,   \Delta t \, \beta_1 
\,+\,   \Delta t^2  \, \beta_2 \,+ \dots \,   + \Delta t^{k} \, \beta_k + \dots \big) \, 
\big( \alpha_1 + \Delta t \, \alpha_2 \,+  \dots    + \Delta t^{\ell} \, \alpha_\ell + \dots\big) ^j 
\, W $

\smallskip  \noindent
relatively to the space                      
derivatives. Then with the help of the induction hypothesis 
\moneq   \label{ap1-serie-Kappa}    
\Big(  E \,+\,  \sum_{m=1}^{\infty}   \Delta t^{m} \, \beta_m  \Big) \, 
\Big(  \sum_{p=1}^{\infty}   \Delta t^{p-1} \, \alpha_p  \Big) ^j \,\equiv \, 
 \sum_{\ell=0}^{\infty}   \Delta t^{\ell} \,\, K^j_{j+\ell} \,  
\, , \quad j \geq 0 \,  ,   
\monend
we identify the two expressions of the coefficient of $  \Delta t^{k+1} $ 
issued from the equation  (\ref{ap1-dvpmt-Y}):   
\moneq   \label{ap1-recurrence-beta}    
S \, \beta_{k+1} \,=\, C_{k+1}  \,+\,  \sum_{j=1}^{k+1} D_j \,\,  \beta_{k+1-j} 
 \,-\,   \sum_{j=1}^{k+1}  {{1}\over{j !}} \,  K^j_{k+1} \, . 
\monend
%

\monitem  The explicitation of the coefficients  $ \, \Gamma^j_{j+\ell} \,$ 
and  $\, K^j_{k+1} \,$  of the matricial formal series   is now easy, due to the relations 
(\ref{ap1-serie-Gamma}) and (\ref{ap1-serie-Kappa}). 
We  specify      
the coefficients $ \,  \Gamma_{j+\ell}^\ell \, $ obtained in the  matricial formal series
 (\ref{ap1-serie-Gamma}). For $j = 0 $, the power in relation 
 (\ref{ap1-serie-Gamma}) is the identity. Then 
\moneq   \label{ap1-Gamma-j0} 
 \Gamma^0_0 = {\rm I} \,, \qquad  \Gamma^0_\ell = 0  \,, \quad  \ell \geq 1 \, .
\monend
When $j=1$, the initial series is not changed. Then 
\moneq   \label{ap1-Gamma-j1} 
 \Gamma^1_\ell = \alpha_\ell \,,  \qquad  \ell \geq 1 \, . 
\monend
For $j=2$, we have to compute the square of the initial series, 
paying attention that the matrix operators $ \, \alpha_\ell \,$      
do not commute.
Observe that with the formal Chapman-Enskog method used {\it e.g.}   in \cite{DdH92}, 
non-commutation relations have also to be taken into consideration 
for higher order terms in the case of several conserved moments. 
  We have

\smallskip \noindent  $ \displaystyle \qquad  \qquad 
\Big(  \sum_{\ell=1}^{\infty}   \Delta t^{\ell} \, \alpha_{\ell+1} \Big) \, 
\Big(  \sum_{j=1}^{\infty}   \Delta t^{j} \, \alpha_{j+1} \Big)  \,=\, 
 \sum_{p=0}^{\infty}  \Delta t^{p} \,  \sum_{\ell + j = p }  \alpha_{\ell+1} \, 
  \alpha_{j+1} $ 

\noindent    and we have in particular
\moneq   \label{ap1-Gamma-j2} 
 \Gamma^2_2 \,=\, \alpha_1^2 \,,\quad 
\Gamma^2_3 \,=\, \alpha_1 \, \alpha_2 \,+\, \alpha_2 \,\alpha_1  \,,\quad  
\Gamma^2_4 \,=\, \alpha_1 \, \alpha_3  \,+\, 
 \alpha_2^2 \,+\,  \alpha_3 \,\alpha_1 \, . 
\monend
In the general case, we have 
\moneqstar  
\Big(  \sum_{\ell=0}^{\infty}   \Delta t^{\ell} \, \alpha_{\ell+1} \Big) ^j \,= \, 
 \sum_{\ell=0}^{\infty}    \,  \Delta t^{p} \,  \sum_{\ell_1 + \dots +\ell_j = p}  
\alpha_{\ell_1+1} \, \dots \, \alpha_{\ell_j+1}   
\monendstar 

\noindent 
and in consequence 
\moneq   \label{ap1-Gamma-Gamma-p} 
\Gamma^j_{p+j} \,=\,       \sum_{\ell_1 + \dots +\ell_j = p}  
\alpha_{\ell_1+1} \, \dots \, \alpha_{\ell_j+1}    \,.
\monend
We have in particular for $j=3$ and $j=4$:
\moneq   \label{ap1-Gamma-j3-4} 
\Gamma^3_3 \,=\, \alpha_1^3 \,, \qquad 
\Gamma^3_4 \,=\, \alpha_1^2 \, \alpha_2 \,+\, 
\alpha_1 \, \alpha_2 \,  \alpha_1 \,+ \,  \alpha_2 \, \alpha_1^2  \,, \qquad 
\Gamma^4_4 \,=\, \alpha_1^4 \, . 
\monend
%
For the explicitation of the coefficients  $\, K^j_{k+1} , \, $  we can replace the power 
of the formal series of the relation  (\ref{ap1-serie-Gamma}) 
in the relation (\ref{ap1-serie-Kappa}).  
We obtain, with the notation $ \, \beta_0 \equiv E$,  
\moneqstar
\Big(  \sum_{m=0}^{\infty}   \Delta t^{m} \, \beta_m \Big) \, 
\Big(  \sum_{\ell=0}^{\infty}   \Delta t^{\ell} \,\, \Gamma^j_{j+\ell} \Big) \,\equiv \, 
 \sum_{p=0}^{\infty}   \Delta t^{p} \,\, K_{j+p}^j   
\monendstar 
then we have by induction 
\moneq   \label{ap1-recurr-Kappa}  
  K_{j+p}^j \,=\,  \sum_{m+\ell=p}  \beta_m \,\, \Gamma^{j}_{j+\ell} \, . 
\monend
For $j=0$, we deduce 
\moneq   \label{ap1-Kappa-j0} 
  K^0_0 = E \,, \qquad K^0_p= 0  \,,  \quad  p \geq 1 \,  
\monend
and for $j=1$, we have a simple product of two formal series:
\moneq   \label{ap1-Kappa-j1} 
 K^1_p  = E \, \alpha_p \,+\,  \beta_1  \, \alpha_{p-1} \,+\, \dots 
\,+\,  \beta_{p-1}  \, \alpha_1  \,,    \quad  p \geq 1 \, . 
\monend
We specify     
 some particular values of the coefficients 
$ \,   K_{j+p}^j \,$ when $j=2$, $j=3$ and for $j=4$:
\moneq   \label{ap1-Kappa-j234} 
\left\{ \begin{array} [c]{l}   \displaystyle 
K^2_2 \,=\, E \, \Gamma^2_2 \,, \quad 
 K^2_3 \,=\, E \, \Gamma^2_3   \,+\, \beta_1 \, \Gamma^2_2  \,, \quad 
 K^2_4 \,=\, E \, \Gamma^2_4  
\,+\, \beta_1 \, \Gamma^2_3 \,+\, \beta_2 \, \Gamma^2_2 \,, 
 \\ \displaystyle   \vspace{-.5cm}  ~  \\ \displaystyle
K^3_3 \,=\, E \, \Gamma^3_3 \,, \quad 
 K^3_4 \,=\, E \, \Gamma^3_4   \,+\, \beta_1 \, \Gamma^3_3  \,, \quad 
 K^4_4  \,=\, E \, \Gamma^4_4 \, .  
\end{array} \right. \monend  
%


\monitem   It is now possible to make explicit      
 up to fourth order to fix the ideas 
the matricial coefficients of the
 expansion  (\ref{ap1-developpement-Y}) of the nonconserved moments 
and of the associated   partial differential equation (\ref{ap1-edp-W}).  
We have,   following the natural order of the algorithm:
\moneq   \label{ap1-algo-ordre-4}  
\left\{ \begin{array} [c]{rcl}   \displaystyle  
\beta_0  & = &  E 
 \\ \displaystyle   \vspace{-.5cm}  ~  \\ \displaystyle 
\alpha_1 & = &  A_1 \,+\, B_1 \, E 
 \\ \displaystyle   \vspace{-.5cm}  ~  \\ \displaystyle 
 \beta_1 &= &  S^{-1} \,  \big( C_1 \,+\, D_1 \, E  \,-\, K^1_1  \big) 
 \\ \displaystyle   \vspace{-.5cm}  ~  \\ \displaystyle
 \alpha_2  & = &    A_2 +  B_2 \, E \,+\,  B_1 \,  \beta_1 \,-\, {1\over2} \,  \Gamma_2^2
 \\ \displaystyle   \vspace{-.5cm}  ~  \\ \displaystyle 
\beta_2  & = &  S^{-1} \, \Big[ \, 
  C_2 \,+\, D_2 \, E \,+\, D_1 \, \beta_1 \,-\,    K^1_2 \,-\, {{1}\over{2}} \, K_2^2  \, \Big]
 \\ \displaystyle   \vspace{-.5cm}  ~  \\ \displaystyle  
\alpha_3  & =  &    A_3 +  B_1 \,  \beta_2 \,+\, B_2 \,  \beta_1 \,+\, B_3 \, E  
\, - \,  {1\over2} \, \Gamma^2_3 \, -  {1\over6} \, \Gamma_3^3 
 \\ \displaystyle   \vspace{-.5cm}  ~  \\ \displaystyle   
\beta_3  & = &  S^{-1} \, \Big[ \, 
   C_3 \,+\, D_1  \,\beta_2 \, + \, D_2  \,\beta_1 \, + \, D_3  \, E 
\,-\, K^1_3 \, - \,  {1\over2} \, K^2_3    \, - \,  {1\over6} \, K^3_3    \, \Big]
 \\ \displaystyle   \vspace{-.5cm}  ~  \\ \displaystyle  
\alpha_4  & = &    A_4 +  B_1 \,  \beta_3 \,+\, B_2 \,  \beta_2  \,+\, B_3 \,  \beta_1  \,+\, B_4 \, E  
\, - \,  {1\over2} \, \Gamma^2_4 \, - \,  {1\over6} \, \Gamma^3_4 \, 
- \,  {1\over24} \, \Gamma_4^4  \, .  
\end{array} \right. \monend  
Observe that with the  explicit relations (\ref{ap1-algo-ordre-4}), the  computer time 
for deriving  formally the equivalent partial equation like  (\ref{ap1-hyp-recur}) 
at fourth order of accuracy   
has been reduced by three orders of magnitude~(!) 
in  comparison  with the algorithm presented in the contribution~\cite{DL09}.


\bigskip \bigskip     
\noindent {\bf \large  Appendix 2. A specific  algebraic coefficient   } 

\monitem  
With   H\'enon's  coefficients $ \, \sigma_j \, $ defined according to (\ref{sigmas}), 
a numbering of the D3Q27 moments proposed in (\ref{d3q27-moments}),
the equilibrium  of the energy   (\ref{d2q9-cns-ordre-1}) parametrized by $\, \alpha $, 
and   the parameter $ \, \psi \,$ introduced in (\ref{d3q27-psi-ordre-4}), 
the coefficient $\, \zeta_4 \, $
for the fourth order term in (\ref{edp-ordre-4}) can be evaluated according to:  
%
%

\newpage 

  \moneq   \label{d2q27-zeta4-ordre-4} 
    \left\{   \begin{array}{l}  \displaystyle  
\zeta_4 \, = \,  {{1}\over{108}} \, {{\sigma_5 \, \, \lambda \,\,  \Delta x^3}\over
{ (4 \, \psi-7) \, (3 \, \psi^2-11 \, \psi +\, 14)
\,  ( 14   +\, 23 \, \psi - 22 \, \psi^2  +\,  3 \, \psi^3)^3}} \, \, N_4 
\\      \displaystyle    
N_4 \, =  - 526848 -\, 13105344 \, \sigma_5^2 +\, 56334931 \, 
    \psi^6 \, \alpha -\, 3413088 \, \psi^2 
\\                                           \displaystyle    \qquad  
+\, 29925576 \, \psi^3
 +\, 44310000 \, \sigma_5^2 \,  \psi^3 \, \alpha^2
 +\,  2458624 \, \alpha^2 
-\, 7776 \, \sigma_5^2 \, \psi^{12} 
\\                                           \displaystyle    \qquad   
+\, 116153808 \, \sigma_5^2 \, \psi^7  
+\, 16213680 \,    \sigma_5^2 \, \psi^9 
 -\, 56871552 \, \sigma_5^2 \, \psi^8 
-\, 2696976 \, \sigma_5^2 \, \psi^{10} 
\\                                           \displaystyle    \qquad  
+\, 803992 \, \psi \,      \alpha 
-\, 16250948 \, \psi^2 \, \alpha 
+\, 15057742 \, \psi^3 \, \alpha 
 +\, 236520 \, \sigma_5^2 \, \psi^{11} 
\\                                           \displaystyle    \qquad  
-\,    47554008 \, \sigma_5^2 \, \psi^5 \, \alpha^2  
+\, 414648 \, \sigma_5^2 \, \psi^9 \, \alpha^2 
+\, 5924856 \,  \sigma_5^2 \, \psi^7 \, \alpha^2  
\\                                           \displaystyle    \qquad 
-\, 3520104 \, \sigma_5^2 \, \psi^8 \, \alpha^2 
 +\, 7776 \, \sigma_5^2 \, \psi^{12} \,  \alpha^2 
-\, 73224 \, \sigma_5^2 \, \psi^{11} \, \alpha^2 
-\, 1805156 \, \psi^8 \, \alpha^2  
\\                                           \displaystyle    \qquad 
-\, 1316084 \, \psi^7 \,     \alpha^2 
+\, 13802956 \, \psi^6 \, \alpha^2
 -\, 29063324 \, \psi^5 \, \alpha^2 
+\, 25708132 \, \psi^4 \, \alpha^2 
\\                                           \displaystyle    \qquad 
-\,   1230152 \, \psi^3 \, \alpha^2 
+\, 3742816 \, \alpha^2 \, \psi 
+\, 27756 \, \psi^{11} \, \alpha^2 
-\, 3429 \, \psi^{11} \,    \alpha 
\\                                           \displaystyle    \qquad
-\, 187851264 \, \sigma_5^2 \, \alpha^2 \, \psi^2 
+\, 198524928 \, \sigma_5^2 \, \alpha^2 \, \psi  
+\, 195048  \, \sigma_5^2 \, \psi^{10} \, \alpha^2 
\\                                           \displaystyle    \qquad
-\, 12827088 \, \psi^4 
+\, 100016448 \, \sigma_5^2 \, \psi 
-\, 10762392 \,     \sigma_5^2 \, \psi^3 
\\                                           \displaystyle    \qquad
 -\, 287184 \, \sigma_5^2 \, \psi^5 
-\, 117365232 \, \sigma_5^2 \, \psi^6 
+\, 102921792 \,     \sigma_5^2 \, \psi^4 
-\, 131926368 \, \sigma_5^2 \, \psi^2 
\\                                           \displaystyle    \qquad
-\, 6082272 \, \psi  
+\, 22678777 \, \psi^4 \, \alpha 
-\,   58798343 \, \psi^5 \, \alpha 
+\, 316283520 \, \sigma_5^2 \, \psi \, \alpha 
\\                                           \displaystyle    \qquad
-\, 421440000 \, \sigma_5^2 \, \psi^2 \,     \alpha 
+\, 148286280 \, \sigma_5^2 \, \psi^3 \, \alpha 
+\, 2458624 \, \alpha -\, 1296 \, \psi^{12} \, \alpha^2  
\\                                           \displaystyle    \qquad 
-\, 324 \,     \psi^{12} \, \alpha 
-\, 12657680 \, \psi^2 \, \alpha^2  
+\, 169965 \, \psi^{10} \, \alpha 
-\, 1834629 \, \psi^9 \, \alpha 
\\                                           \displaystyle    \qquad 
+\,   9787591 \, \psi^8 \, \alpha 
-\, 30298973 \, \psi^7 \, \alpha 
-\, 235980 \, \psi^{10} \, \alpha^2 
+\, 989292 \, \psi^9 \,    \alpha^2 
\\                                           \displaystyle    \qquad 
+\, 55400808 \, \sigma_5^2 \, \psi^4 \, \alpha^2  
-\, 3888 \, \sigma_5^2 \, \psi^{12} \, \alpha
 +\, 223236 \,    \sigma_5^2 \, \psi^{11} \, \alpha  
\\                                           \displaystyle    \qquad 
-\, 2725812 \, \sigma_5^2 \, \psi^{10} \, \alpha   
+\, 15619428 \, \sigma_5^2 \, \psi^9 \, \alpha  
-\, 48491916 \, \sigma_5^2 \, \psi^8 \, \alpha 
\\                                           \displaystyle    \qquad  
+\, 75338436 \, \sigma_5^2 \, \psi^7 \, \alpha  
+\, 8771448 \,   \sigma_5^2 \, \psi^6 \, \alpha^2  
-\, 78989568 \, \sigma_5^2 \, \alpha 
\\                                           \displaystyle    \qquad 
-\, 6400920 \, \psi^9 -\, 56568924 \, \psi^7 
+\,      24275088 \, \psi^8 +\, 3240 \, \psi^{12} 
-\, 88452 \, \psi^{11}  
\\                                           \displaystyle    \qquad 
-\, 46884384 \, \psi^5 +\, 76942908 \, \psi^6 
+\,       151481100 \, \sigma_5^2 \, \psi^4 \, \alpha 
-\, 12989436 \, \sigma_5^2 \, \psi^6 \, \alpha 
\\                                           \displaystyle    \qquad 
-\, 141331668 \,   \sigma_5^2 \, \psi^5 \, \alpha 
+\, 1015308 \, \psi^{10} 
-\, 77070336 \, \sigma_5^2 \,   \alpha^2 \, . 
\end{array}\right. \monend

\bigskip  

  \end{document}